
%
%
%
%
%
%
\magnification=\magstephalf      
%
%
\vsize=7.5truein                 
\hsize=5.2truein                 
\newskip\stdskip                 
\stdskip=6pt plus3pt minus3pt    
\medskipamount=\stdskip          
\parindent=0pt                   
\parskip=\stdskip                
\abovedisplayskip=\stdskip       
\belowdisplayskip=\stdskip       
\mathsurround=0.75pt             
\overfullrule=0pt                
%
%
\def\ppar{\par\goodbreak\vskip 8pt plus 4pt minus 4pt}     
%
%
\def\stdspace{\hskip 0.75em plus 0.15em\ignorespaces}
\let\qua\stdspace 
%
%
%
%
%
%
%
\def\hexnumber#1{\ifcase#1 0\or 1\or 2\or 3\or 4\or 5\or 6\or 7\or 8\or
 9\or A\or B\or C\or D\or E\or F\fi}
%
%
\font\thirtnmsa=msam10 scaled 1315    
\font\tenmsa=msam10          \font\ninemsa=msam9
\font\sevenmsa=msam7         \font\sixmsa=msam6
\font\fivemsa=msam5
%
%
\newfam\msafam                  \textfont\msafam=\tenmsa
\scriptfont\msafam=\sevenmsa    \scriptscriptfont\msafam=\fivemsa
\edef\hexa{\hexnumber\msafam}        
\def\msa{\fam\msafam\tenmsa}         
%
%
\font\thirtnmsb=msbm10 scaled 1315   
\font\tenmsb=msbm10      \font\ninemsb=msbm9
\font\sevenmsb=msbm7     \font\sixmsb=msbm6
\font\fivemsb=msbm5
%
\newfam\msbfam                   \textfont\msbfam=\tenmsb       
\scriptfont\msbfam=\sevenmsb     \scriptscriptfont\msbfam=\fivemsb
\edef\hexb{\hexnumber\msbfam}    
\def\msb{\fam\msbfam\tenmsb}     
%
%
\font\thirtneufm=eufm10 scaled 1315   
\font\teneufm=eufm10                 \font\nineeufm=eufm9
\font\seveneufm=eufm7                \font\sixeufm=eufm6
\font\fiveeufm=eufm5
%
\newfam\eufmfam                    \textfont\eufmfam=\teneufm
\scriptfont\eufmfam=\seveneufm     \scriptscriptfont\eufmfam=\fiveeufm
\edef\hexf{\hexnumber\eufmfam}      
\def\frak{\fam\eufmfam\teneufm}     
%
%
%
\font\thirtnrm=cmr10 scaled 1315    
\font\ninerm=cmr9                   \font\sixrm=cmr6   
%
\font\thirtni=cmmi10 scaled 1315    
\font\ninei=cmmi9                   \font\sixi=cmmi6  
%
\font\thirtnsy=cmsy10 scaled 1315   
\font\ninesy=cmsy9                  \font\sixsy=cmsy6  
%
\font\thirtnbf=cmbx10 scaled 1315   
\font\ninebf=cmbx9                  \font\sixbf=cmbx6  
%
%
\font\thirtnex=cmex10 scaled 1315   
\font\nineex=cmex9                  
%
%
\font\thirtnit=cmti10 scaled 1315  
\font\nineit=cmti9                  
%
\font\thirtnsl=cmsl10 scaled 1315  
\font\ninesl=cmsl9                  
%
\font\thirtntt=cmtt10 scaled 1315  
\font\ninett=cmtt9                  
%
%
%
%
\def\small{%
%
%
\textfont0=\ninerm \scriptfont0=\sixrm \scriptscriptfont0=\fiverm
\def\rm{\fam0\ninerm}
%
%
\textfont1=\ninei \scriptfont1=\sixi \scriptscriptfont1=\fivei
%
%
\textfont2=\ninesy \scriptfont2=\sixsy \scriptscriptfont2=\fivesy
%
%
\textfont3=\nineex \scriptfont3=\nineex \scriptscriptfont3=\nineex
%
%
\textfont\bffam=\ninebf \scriptfont\bffam=\sixbf
\scriptscriptfont\bffam=\fivebf \def\bf{\fam\bffam\ninebf}%
%
%
\textfont\itfam=\nineit \def\it{\fam\itfam\nineit}%
\textfont\slfam=\ninesl \def\sl{\fam\slfam\ninesl}%
\textfont\ttfam=\ninett \def\tt{\fam\ttfam\ninett}%
%
%
%
\textfont\msafam=\ninemsa \scriptfont\msafam=\sixmsa
\scriptscriptfont\msafam=\fivemsa \def\msa{\fam\msafam\ninemsa}%
%
%
\textfont\msbfam=\ninemsb \scriptfont\msbfam=\sixmsb
\scriptscriptfont\msbfam=\fivemsb \def\msb{\fam\msbfam\ninemsb}%
%
%
\textfont\eufmfam=\nineeufm  \scriptfont\eufmfam=\sixeufm
\scriptscriptfont\eufmfam=\fiveeufm \def\frak{\fam\eufmfam\nineeufm}%
%
%
%
\normalbaselineskip=11pt%
\setbox\strutbox=\hbox{\vrule height8pt depth3pt width0pt}%
%
%
\normalbaselines\rm
%
%
\stdskip=4pt plus2pt minus2pt    
\medskipamount=\stdskip          
\parskip=\stdskip                
\abovedisplayskip=\stdskip       
\belowdisplayskip=\stdskip       
\def\ppar{\par\goodbreak\vskip 6pt plus 3pt minus 3pt}%
%
%
\def\section##1{\global\advance\sectionnumber by 1
\vskip-\lastskip\penalty-800\vskip 20pt plus10pt minus5pt 
\egroup{\bf\number\sectionnumber\quad##1}\bgroup\small         
\vskip 6pt plus3pt minus3pt
\nobreak\resultnumber=1}
}    
%
\def\beginsmall{\bgroup\small}
\let\endsmall\egroup
%
%
%
%
\def\large{%
\textfont0=\thirtnrm \scriptfont0=\ninerm \scriptscriptfont0=\sevenrm
\def\rm{\fam0\thirtnrm}%
\textfont1=\thirtni \scriptfont1=\ninei \scriptscriptfont1=\seveni
\textfont2=\thirtnsy \scriptfont2=\ninesy \scriptscriptfont2=\sevensy
\textfont3=\thirtnex \scriptfont3=\thirtnex \scriptscriptfont3=\thirtnex
\textfont\bffam=\thirtnbf \scriptfont\bffam=\ninebf
\scriptscriptfont\bffam=\sevenbf \def\bf{\fam\bffam\thirtnbf}%
\textfont\itfam=\thirtnit \def\it{\fam\itfam\thirtnit}%
\textfont\slfam=\thirtnsl \def\sl{\fam\slfam\thirtnsl}%
\textfont\ttfam=\thirtntt \def\tt{\fam\ttfam\thirtntt}%
\textfont\msafam=\thirtnmsa \scriptfont\msafam=\ninemsa
\scriptscriptfont\msafam=\sevenmsa \def\msa{\fam\msafam\thirtnmsa}%
\textfont\msbfam=\thirtnmsb \scriptfont\msbfam=\ninemsb
\scriptscriptfont\msbfam=\sevenmsb \def\msb{\fam\msbfam\thirtnmsb}%
\textfont\eufmfam=\thirtneufm  \scriptfont\eufmfam=\nineeufm
\scriptscriptfont\eufmfam=\seveneufm \def\frak{\fam\eufmfam\teneufm}%
\normalbaselineskip=16pt%
\setbox\strutbox=\hbox{\vrule height11.5pt depth4.5pt width0pt}%
\normalbaselines\rm}%
\let\Large\large   
%

%

%
\mathchardef\plussquare="0\hexa01
\mathchardef\nge="3\hexb0B
\mathchardef\maltesecross="0\hexa7A
\mathchardef\del="0\hexf01
%
%
%
%
\font\sc=cmcsc10
%
%
%
%
\def\sqr#1#2{{\vcenter{\vbox{\hrule  height.#2truept
	\hbox{\vrule width.#2truept height#1truept 
	\kern#1truept \vrule width.#2truept}
	\hrule height.#2truept}}}}
\def\sq{\sqr55}    
%
%
%
%
\newcount\sectionnumber            
\newcount\resultnumber             
\sectionnumber=0\resultnumber=1    
%
%
%
\def\section#1{\global\advance\sectionnumber by 1
\xdef\nextkey{\number\sectionnumber}
\vskip-\lastskip\penalty-800\vskip 20pt plus10pt minus5pt 
{\large\bf\number\sectionnumber\quad#1}         
\vskip 8pt plus4pt minus4pt
\nobreak\resultnumber=1}      
%
%
%
%
%
\def\sh#1{\vskip-\lastskip\ppar{\bf #1}\par\nobreak\medskip}         
%
%
%
%

%
\def\proc#1{\xdef\nextkey{\number\sectionnumber.\number\resultnumber}%
\vskip-\lastskip\ppar\bf%
\noindent#1\ \number\sectionnumber.\number\resultnumber
\stdspace\sl\global\advance\resultnumber by 1\ignorespaces}
\def\endproc{\rm\ppar} 
%
%
\def\prf{\vskip-\lastskip\ppar\noindent{\bf Proof}%
\stdspace\rm}                            
\def\qed{\hfill$\sq$\par\goodbreak\rm}   
\let\endproof\endprf              
%
%
%
%
%
%
%
%
\def\proclaim#1{\vskip-\lastskip\ppar\bf%
\noindent#1\stdspace\sl\ignorespaces} 

%
%
%
%
\def\rk#1{\vskip-\lastskip\ppar{\bf #1}\stdspace\ignorespaces}                

%
%
%
%
%
%
\def\label{\xdef\nextkey{\number\sectionnumber.\number\resultnumber}%
\number\sectionnumber.\number\resultnumber
\global\advance\resultnumber by 1}
%
%
%
%
%
%
%
%
%
%
%
%
%
%
%
%
\newcount\refnumber              
\refnumber=1                     
\long\def\reflist#1\endreflist{%
\long\def\thereflist{#1}{\def\refkey##1##2\par{\xdef##1{\number\refnumber}%
\global\advance\refnumber by 1}%
\def\key##1##2\par{\expandafter\xdef%
\csname##1\endcsname{\number\refnumber}%
\global\advance\refnumber by 1}#1\par}}
\long\def\references{%
\penalty-800\vskip-\lastskip\vskip 15pt plus10pt minus5pt 
{\large\bf References}\ppar 
{\leftskip=25pt\frenchspacing    
\small\parskip=3pt plus2pt       
\def\refkey##1##2\par{\noindent  
\llap{[##1]\stdspace}\ignorespaces##2\par}         
\def\key##1##2\par{\noindent  
\llap{[\ref{##1}]\stdspace}\ignorespaces##2\par}  
\def\,{\thinspace}\thereflist\par}}
%
%
%
\newcount\footnotenumber         
\footnotenumber=1                
\def\fnote#1{\xdef\nextkey{\number\footnotenumber}%
{\small\ifnum\footnotenumber>9\parindent=14pt%
\else\parindent=10pt\fi\footnote{$^{\number\footnotenumber}$}%
{\hglue-5pt#1}\global\advance\footnotenumber by 1}}
%
%
%
%
%
%
%
\newcount\figurenumber          
\figurenumber=1                 
\def\caption#1{\xdef\nextkey{\number\figurenumber}%
\cl{\small Figure \number\figurenumber: #1}%
\global\advance\figurenumber by 1}
\def\figurelabel{\xdef\nextkey{\number\figurenumber}%
\cl{\small Figure \number\figurenumber}%
\global\advance\figurenumber by 1}
\long\def\figure#1\endfigure{{\xdef\nextkey{\number\figurenumber}%
\let\captiontext\relax\def\caption##1{\xdef\captiontext{##1}}%
\midinsert\cl{\ignorespaces#1\unskip\unskip\unskip\unskip}\vglue6pt\cl{\small 
Figure \number\figurenumber\ifx\captiontext\relax\else: \captiontext
\fi}\endinsert\global\advance\figurenumber by 1}}
%
%
%
%
%
%
%
\def\nextkey{??}   
%
\def\key#1{\expandafter\xdef\csname #1\endcsname{\nextkey}}
\def\ref#1{\expandafter\ifx\csname #1\endcsname\relax
\immediate\write16{Reference {#1} undefined}??\else
\csname #1\endcsname\fi}
%
%
%
%
%
%
%
\newread\gtinfile
\newwrite\gtreffile
\def\useforwardrefs{
\openin\gtinfile\jobname.ref
\ifeof\gtinfile
\closein\gtinfile
\immediate\write16{No file \jobname.ref}
\else
\closein\gtinfile
\input \jobname.ref
\fi
\immediate\openout\gtreffile \jobname.ref
%
%
\def\key##1{{\def\\{\noexpand}%
\expandafter\xdef\csname ##1\endcsname{\nextkey}%
\immediate\write\gtreffile{\\\expandafter\\\def\\\csname ##1\\\endcsname%
{\nextkey}}}}
%
%
\long\def\reflist##1\endreflist{%
\long\def\thereflist{##1}{\def\refkey####1####2\par{\xdef####1{%
\number\refnumber}{\def\\{\noexpand}\immediate\write\gtreffile
{\\\def\\####1{\number\refnumber}}}\global\advance\refnumber by 1}%
\def\key####1####2\par{\expandafter\xdef%
\csname####1\endcsname{\number\refnumber}%
{\def\\{\noexpand}\immediate\write\gtreffile
{\\\expandafter\\\def\\\csname ####1\\\endcsname{\number\refnumber}}}
\global\advance\refnumber by 1}##1\par}}
\long\def\biblio##1\endbiblio{\reflist##1\endreflist\references}%
%
%
\def\numkey##1{{\def\\{\noexpand}%
\xdef##1{\number\sectionnumber.\number\resultnumber}
\immediate\write\gtreffile{\\\def\\##1%
{\number\sectionnumber.\number\resultnumber}}}}
\def\seckey##1{{\def\\{\noexpand}\xdef##1{\number\sectionnumber}
\immediate\write\gtreffile{\\\def\\##1{\number\sectionnumber}}}}
\def\figkey##1{\xdef##1{\number\figurenumber}%
{\def\\{\noexpand}\immediate\write\gtreffile%
{\\\def\\##1{\number\figurenumber}}}
\number\figurenumber\global\advance\figurenumber by 1}
}   
%
%
%
%
\def\figkey#1{\xdef#1{\number\figurenumber}%
\number\figurenumber\global\advance\figurenumber by 1}
\def\fig#1#2\endfig{%
\midinsert\cl{#2}\vglue6pt\cl{\small Figure #1}\endinsert}
\def\newfig{\number\figurenumber\global\advance\figurenumber by 1}
\def\numkey#1{\xdef#1{\number\sectionnumber.\number\resultnumber}}
\def\seckey#1{\xdef#1{\number\sectionnumber}}
%
%
%
%
%
%
%
%
%
\def\verb{\catcode`\"=\active}       
\def\brev{\catcode`\"=12}            
\brev                                
\verb                                
{\obeyspaces\gdef {\ }}              
{\catcode`\`=\active\gdef`{\relax\lq}}
\def"{%
\begingroup\baselineskip=12pt\def\par{\leavevmode\endgraf}%
\tt\obeylines\obeyspaces\parskip=0pt\parindent=0pt%
\catcode`\$=12\catcode`\&=12\catcode`\^=12\catcode`\#=12%
\catcode`\_=12\catcode`\~=12%
\catcode`\{=12\catcode`\}=12\catcode`\%=12\catcode`\\=12%
\catcode`\`=\active\let"\endgroup}
\brev      
%
%
%
%
%
%
\def\items{\par\leftskip = 25pt}           
\def\enditems{\par\leftskip = 0pt}         
\def\item#1{\par\leavevmode\llap{#1\stdspace}%
\ignorespaces}                             
%
%

%
%
\def\co{\colon\thinspace}    
\def\np{\vfil\eject}         
\def\nl{\hfil\break}         
\def\cl{\centerline}         
\def\gt{{\mathsurround=0pt\it $\cal G\mskip-2mu$eometry \&\ 
$\cal T\!\!$opology}}        
\def\agt{{\mathsurround=0pt\it$\cal A\mskip-.7mu$lgebraic \&\ 
$\cal G\mskip-2mu$eometric $\cal T\!\!$opology}}  
%
%
%

%
%
%
%
%
\def\title#1{\def\thetitle{#1}}

\def\author#1{\edef\previousauthors{\theauthors}
 \ifx\theauthors\relax\def\theauthors{#1}\else
 \def\theauthors{\previousauthors\par#1}\fi}

%
\def\address#1{\edef\previousaddresses{\theaddress}
 \ifx\theaddress\relax\def\theaddress{#1}\else
 \def\theaddress{\previousaddresses\par\vskip 2pt\par#1}\fi}
\def\secondaddress#1{\edef\previousaddresses{\theaddress}
 \ifx\theaddress\relax\def\theaddress{#1}\else
 \def\theaddress{\previousaddresses\par{\rm and}\par#1}\fi}   

\def\email#1{\edef\previousemails{\theemail}
 \ifx\theemail\relax\def\theemail{#1}\else
 \def\theemail{\previousemails\hskip 0.75em\relax#1}\fi}
\def\secondemail#1{\edef\previousemails{\theemail}
 \ifx\theemail\relax\def\theemail{#1}\else
 \def\theemail{\previousemails\hskip 0.75em{\rm and}\hskip 0.75em
 \relax#1}\fi}
\def\url#1{\edef\previousurls{\theurl}
 \ifx\theurl\relax\def\theurl{#1}\else
 \def\theurl{\previousurls\hskip 0.75em\relax#1}\fi}
\def\secondurl#1{\edef\previousurls{\theurl}
 \ifx\theurl\relax\def\theurl{#1}\else
 \def\theurl{\previousurls\hskip 0.75em{\rm and}\hskip 0.75em
 \relax#1}\fi}
\long\def\abstract#1\endabstract{\long\def\theabstract{#1}}
\def\primaryclass#1{\def\theprimaryclass{#1}}
\def\secondaryclass#1{\def\thesecondaryclass{#1}}
\def\keywords#1{\def\thekeywords{#1}}
%
%
\let\\\par\let\thetitle\relax\let\theshorttitle\relax
\let\theauthors\relax\let\theshortauthors\relax
\let\theaddress\relax\let\theshortaddress\relax
\let\theemail\relax\let\theurl\relax
\let\theabstract\relax\let\theprimaryclass\relax
\let\thesecondaryclass\relax\let\thekeywords\relax
%
%
%
%
\long\def\maketitlepage{    

\vglue 0.2truein   

%
{\parskip=0pt\leftskip 0pt plus 1fil\def\\{\par\smallskip}{\large
\bf\thetitle}\par\medskip}   

\vglue 0.15truein 

%
{\parskip=0pt\leftskip 0pt plus 1fil\def\\{\par}{\sc\theauthors}
\par\medskip}%
 
\vglue 0.1truein 

%
{\small\parskip=0pt
{\leftskip 0pt plus 1fil\def\\{\par}{\sl\theaddress}\par}
\ifx\theemail\relax\else  
\vglue 5pt \def\\{\stdspace{\rm and}\stdspace} 
\cl{Email:\stdspace\tt\theemail}\fi
\ifx\theurl\relax\else    
\vglue 5pt \def\\{\stdspace{\rm and}\stdspace} 
\cl{URL:\stdspace\tt\theurl}\fi\par}

\vglue 7pt 

{\bf Abstract}

\vglue 5pt

\theabstract

\vglue 7pt 

{\bf AMS Classification numbers}\quad Primary:\quad \theprimaryclass\par

Secondary:\quad \thesecondaryclass

\vglue 5pt 

{\bf Keywords:}\quad \thekeywords

\np  

}    
%
%
\long\def\makeshorttitle{    


%
{\parskip=0pt\leftskip 0pt plus 1fil\def\\{\par\smallskip}{\large
\bf\thetitle}\par\medskip}   

\vglue 0.05truein 

%
{\parskip=0pt\leftskip 0pt plus 1fil\def\\{\par}{\sc\theauthors}
\par\medskip}%
 
\vglue 0.03truein 

%
{\small\parskip=0pt
{\leftskip 0pt plus 1fil\def\\{\par}{\sl\ifx\theshortaddress\relax
\theaddress\else\theshortaddress\fi}\par}
\ifx\theemail\relax\else  
\vglue 5pt \def\\{\stdspace{\rm and}\stdspace} 
\cl{Email:\stdspace\tt\theemail}\fi
\ifx\theurl\relax\else    
\vglue 5pt \def\\{\stdspace{\rm and}\stdspace} 
\cl{URL:\stdspace\tt\theurl}\fi\par}

\vglue 10pt 


{\small\leftskip 25pt\rightskip 25pt{\bf Abstract}\stdspace\theabstract

{\bf AMS Classification}\stdspace\theprimaryclass
\ifx\thesecondaryclass\relax\else; \thesecondaryclass\fi\par
{\bf Keywords}\stdspace \thekeywords\par}
\vglue 7pt
}    
\let\maketitle\makeshorttitle        
%
%

\def\volumenumber#1{\def\thevolumenumber{#1}}
\def\volumeyear#1{\def\thevolumeyear{#1}}
\def\pagenumbers#1#2{\def\startpage{#1}\def\finishpage{#2}}
\def\published#1{\def\publishdate{#1}}
\def\received#1{\def\receiveddate{#1}}
\def\revised#1{\def\reviseddate{#1}}
\let\reviseddate\relax
\volumenumber{X}
\volumeyear{20XX}
\pagenumbers{1}{XXX}
\published{XX Xxxember 20XX}

\long\def\makeagttitle{   
\agt\hfill      
\hbox to 60truept{\vbox to 0pt{\vglue -14truept{\bf [Logo here]}\vss}\hss}
\break
{\small Volume \thevolumenumber\ (\thevolumeyear)
\startpage--\finishpage\nl
Published: \publishdate}

\vglue .2truein

{\parskip=0pt\leftskip 0pt plus 1fil\def\\{\par\smallskip}{\large
\bf\thetitle}\par\medskip}   
\vglue 0.05truein 

%
{\parskip=0pt\leftskip 0pt plus 1fil\def\\{\par}{\sc\theauthors}
\par\medskip}%
 
\vglue 0.03truein 


{\small\leftskip 25truept\rightskip 25truept{\bf Abstract}\stdspace\theabstract

{\bf AMS Classification}\stdspace\theprimaryclass
\ifx\thesecondaryclass\relax\else; \thesecondaryclass\fi\par
{\bf Keywords}\stdspace \thekeywords\par}\vglue 7truept

}   


\def\Addresses{\bigskip
{\small \parskip 0pt \leftskip 0pt \rightskip 0pt plus 1fil \def\\{\par}
\sl\theaddress\par\medskip \rm Email:\stdspace\tt\theemail\par
\ifx\theurl\relax\else\smallskip \rm URL:\stdspace\tt\theurl\par\fi}}

\def\agtart{
\hoffset 14truemm
\voffset 31truemm
\font\phead=cmsl9 scaled 950
\font\pnum=cmbx10 scaled 913
\font\pfoot=cmsl9 scaled 950
\headline{\vbox to 0pt{\vskip -4.5mm\line{\small\phead\ifnum
\count0=\startpage ISSN numbers are printed here
\hfill {\pnum\folio}\else\ifodd\count0\def\\{ }%
\ifx\theshorttitle\relax\thetitle\else\theshorttitle\fi\hfill{\pnum\folio}
\else\def\\{ and }{\pnum\folio}\hfill\ifx\theshortauthors\relax\theauthors
\else\theshortauthors\fi\fi\fi}\vss}}
\footline{\vbox to 0pt{\vglue 0mm\line{\small\pfoot\ifnum\count0=\startpage
Copyright declaration is printed here\hfill\else
\agt, Volume \thevolumenumber\ (\thevolumeyear)\hfill\fi}\vss}}
\let\maketitle\makeagttitle\let\makeshorttitle\makeagttitle}


\def\ifplaintex{\expandafter\ifx\csname documentclass\endcsname\relax}


\ifplaintex 
\hoffset 14truemm
\voffset 31truemm
\else
\headsep 23pt
\footskip 35pt
\hoffset -4truemm
\voffset 12.5truemm
\fi

\expandafter\ifx\csname beginpicture\endcsname\relax
\expandafter\ifx\csname documentclass\endcsname\relax
\input pictex \else\font\fiverm=cmr5
\input prepictex \input pictex \input postpictex \fi\fi

\def\gt{{\mathsurround=0pt\it $\cal G\mskip-2mu$eometry \&\ 
$\cal T\!\!$opology}}        

\def\gtp{{\mathsurround=0pt\it $\cal G\mskip-2mu$eometry \&\ 
$\cal T\!\!$opology $\cal P\!$ublications}}  


\def\lognumber#1{\def\thelognumber{#1}}
\def\volumenumber#1{\def\thevolumenumber{#1}}
\def\papernumber#1{\def\thepapernumber{#1}}
\def\volumeyear#1{\def\thevolumeyear{#1}}

\def\pagenumbers#1#2{\def\startpage{#1}\def\finishpage{#2}}
\def\published#1{\def\publishdate{#1}}
\def\proposed#1{\def\theproposer{#1}}
\def\seconded#1{\def\theseconders{#1}}
\def\received#1{\def\receiveddate{#1}}
\def\revised#1{\def\reviseddate{#1}}
\def\accepted#1{\def\accepteddate{#1}}

\def\asciiaddress#1{\def\theasciiaddress{#1}}
\def\asciiemail#1{\def\theasciiemail{#1}}
\long\def\asciiabstract#1{\long\def\theasciiabstract{#1}}
\def\asciikeywords#1{\def\theasciikeywords{#1}}


\let\\\par\let\thelognumber\relax
\let\thevolumenumber\relax\let\thepapernumber\relax
\let\thevolumeyear\relax\let\thesamplenumber\relax\let\startpage\relax
\let\finishpage\relax\let\publishdate\relax\let\receiveddate\relax
\let\reviseddate\relax\let\accepteddate\relax\let\theasciititle\relax
\let\theasciiauthors\relax\let\theasciiaddress\relax
\let\theasciiabstract\relax\let\theasciikeywords\relax
\let\theasciiemail\relax\let\theshortauthors\relax\let\theshorttitle\relax

\long\def\maketitlep{   

\count0=\startpage

\gt\hfill      
\beginpicture
\setcoordinatesystem units <0.33truein, 0.33truein> point at 2.2 0.9
\setplotsymbol ({$\cal G$})
\plotsymbolspacing=9truept
\circulararc 315 degrees from 0 1 center at 0 0
\setplotsymbol ({$\cal T$})
\circulararc 315 degrees from 1 -1 center at 1 0
\endpicture
%
\break
{\small\ifx\thesamplenumber\relax 
Volume \else Sample
\fi\thevolumenumber\ (\thevolumeyear)
\startpage--\finishpage\nl
Published: \publishdate}
\vglue 0.5truein plus 0.4fil minus 0.1truein

{\parskip=0pt\leftskip 0pt plus 1fil\def\\{\par\smallskip}{\ifplaintex\large
\else\Large\fi\bf\thetitle}\par\medskip}   

\vglue 0pt plus 0.1fil 

{\parskip=0pt\leftskip 0pt plus 1fil\def\\{\par}{\sc\theauthors}
\par\medskip}

\vglue 0pt plus 0.1fil 

{\small\parskip=0pt\let\newline\\
{\leftskip 0pt plus 1fil\def\\{\par}{\sl\theaddress}\par}
\expandafter\ifx\theemail\relax    
\relax\else\vglue 5pt plus 0.02fil minus 2pt\def\\{\stdspace{\rm 
and}\stdspace} 
\cl{Email:\stdspace\tt\theemail}\fi
\ifx\theurl\relax                  
\relax\else\vglue 5pt plus 0.02fil minus 2pt\def\\{\stdspace{\rm 
and}\stdspace}
\cl{URL:\stdspace\tt\theurl}\fi\par}

\vglue 7pt plus 0.3fil minus 3pt

{\bf Abstract}
\vglue 5pt plus 0.1fil minus 2pt

\theabstract

\vglue 7pt plus 0.3fil minus 3pt

{\bf AMS Classification numbers}\quad Primary:\quad \theprimaryclass

Secondary:\quad \thesecondaryclass

\vglue 5pt plus 0.3fil minus 2pt

{\bf Keywords}\quad \thekeywords

\vglue 10pt plus 0.5fil minus 5pt

{\small  Proposed: \theproposer\hfill Received: \receiveddate\nl
Seconded: \theseconders\hfill 
\ifx\reviseddate\relax                         
Accepted: \accepteddate                        
\else
Revised: \reviseddate                          
\fi}
\eject
}       

\let\maketitlepage\maketitlep
\let\maketitle\maketitlepage


\font\phead=cmsl9 scaled 950
\font\lhead=cmsl9 scaled 1050
\font\pnum=cmbx10 scaled 913
\font\lnum=cmbx10 
\font\pfoot=cmsl9 scaled 950
\font\lfoot=cmsl9 scaled 1050
\ifplaintex
\headline{\vbox to 0pt{\vskip -4.5mm\line{\small\phead\ifnum
\count0=\startpage ISSN 1364-0380 (on line)
1465-3060 (printed) \hfill {\pnum\folio}\else\ifodd\count0\def\\{ }%
\ifx\theshorttitle\relax\thetitle\else\theshorttitle\fi\hfill{\pnum\folio}
\else\def\\{ and }{\pnum\folio}\hfill\ifx\theshortauthors\relax\theauthors
\else\theshortauthors\fi\fi\fi}\vss}}
\footline{\vbox to 0pt{\vglue 0mm\line{\small\pfoot\ifnum\count0=\startpage
\copyright\ \gtp\hfill\else
\gt, Volume \thevolumenumber\ (\thevolumeyear)\hfill\fi}\vss
}}
\else
\makeatletter
\def\@oddhead{{\small\lhead\ifnum\count0=\startpage ISSN 1364-0380 (on line)
1465-3060 (printed) \hfill {\lnum\number\count0}\else\ifodd\count0
\def\\{ }\ifx\theshorttitle\relax \thetitle \else\theshorttitle\fi\hfill
{\lnum\number\count0}\else\def\\{ and }{\lnum\number\count0}
\hfill\ifx\theshortauthors\relax 
\theauthors\else\theshortauthors\fi\fi\fi}}\def\@evenhead{\@oddhead}
\def\@oddfoot{\small\lfoot\ifnum\count0=\startpage\copyright\ \gtp\hfill\else
\gt, Volume \thevolumenumber\ (\thevolumeyear)\hfill\fi}
\def\@evenfoot{\@oddfoot}
\makeatother
\fi


\newwrite\gtoutfile
\long\gdef\makeheadfile{  
{\def\\{, }\def\s{ }
\immediate\openout\gtoutfile head.xxx
\immediate\write\gtoutfile{To: math@arxiv.org}
\immediate\write\gtoutfile{Subject: put or rep NNNNN:pppp}
\immediate\write\gtoutfile{--text follows this line--}
\immediate\write\gtoutfile{Proxy-for: \ifx\theasciiauthors\relax
\theauthors\else\theasciiauthors\fi\s<\ifx\theasciiemail\relax\theemail\else\theasciiemail\fi>}
\immediate\write\gtoutfile{\noexpand\\}
\immediate\write\gtoutfile{Authors: \ifx\theasciiauthors\relax
\theauthors\else\theasciiauthors\fi}
{\def\\{ }\immediate\write\gtoutfile{Title: \ifx\theasciititle\relax
\thetitle\else\theasciititle\fi}}
\immediate\write\gtoutfile{Subj-class: GT or SG or MG etc}
\immediate\write\gtoutfile{MSC-class: \theprimaryclass\ifx\thesecondaryclass\relax\else, \thesecondaryclass\fi}
\immediate\write\gtoutfile{Journal-ref: Geom. Topol. \thevolumenumber
(\thevolumeyear) \startpage-\finishpage}
\immediate\write\gtoutfile{Comments: Published by Geometry and Topology at}
\immediate\write\gtoutfile{\s\s http://www.maths.warwick.ac.uk/gt/GTVol\thevolumenumber/paper\thepapernumber.abs.html}
\immediate\write\gtoutfile{\noexpand\\}
\immediate\write\gtoutfile{}
\ifx\theasciiabstract\relax
\immediate\write\gtoutfile{\theabstract}\else
\immediate\write\gtoutfile{\theasciiabstract}\fi
\immediate\write\gtoutfile{}
\immediate\write\gtoutfile{\noexpand\\}
\immediate\write\gtoutfile{}
\immediate\closeout\gtoutfile}}  

\def\maketitlepage{\maketitlep\makeheadfile}
\let\maketitle\maketitlepage


\def\ifplaintex{\expandafter\ifx\csname documentclass\endcsname\relax}


\ifplaintex 
\hoffset 14truemm
\voffset 31truemm
\else
\headsep 23pt
\footskip 35pt
\hoffset -4truemm
\voffset 12.5truemm
\fi

\expandafter\ifx\csname beginpicture\endcsname\relax
\expandafter\ifx\csname documentclass\endcsname\relax
\input pictex \else\font\fiverm=cmr5
\input prepictex \input pictex \input postpictex \fi\fi

\def\gt{{\mathsurround=0pt\it $\cal G\mskip-2mu$eometry \&\ 
$\cal T\!\!$opology}}        

\def\gtp{{\mathsurround=0pt\it $\cal G\mskip-2mu$eometry \&\ 
$\cal T\!\!$opology $\cal P\!$ublications}}  


\def\lognumber#1{\def\thelognumber{#1}}
\def\volumenumber#1{\def\thevolumenumber{#1}}
\def\papernumber#1{\def\thepapernumber{#1}}
\def\volumeyear#1{\def\thevolumeyear{#1}}

\def\pagenumbers#1#2{\def\startpage{#1}\def\finishpage{#2}}
\def\published#1{\def\publishdate{#1}}
\def\proposed#1{\def\theproposer{#1}}
\def\seconded#1{\def\theseconders{#1}}
\def\received#1{\def\receiveddate{#1}}
\def\revised#1{\def\reviseddate{#1}}
\def\accepted#1{\def\accepteddate{#1}}

\def\asciiaddress#1{\def\theasciiaddress{#1}}
\def\asciiemail#1{\def\theasciiemail{#1}}
\long\def\asciiabstract#1{\long\def\theasciiabstract{#1}}
\def\asciikeywords#1{\def\theasciikeywords{#1}}


\let\\\par\let\thelognumber\relax
\let\thevolumenumber\relax\let\thepapernumber\relax
\let\thevolumeyear\relax\let\thesamplenumber\relax\let\startpage\relax
\let\finishpage\relax\let\publishdate\relax\let\receiveddate\relax
\let\reviseddate\relax\let\accepteddate\relax\let\theasciititle\relax
\let\theasciiauthors\relax\let\theasciiaddress\relax
\let\theasciiabstract\relax\let\theasciikeywords\relax
\let\theasciiemail\relax\let\theshortauthors\relax\let\theshorttitle\relax

\long\def\maketitlep{   

\count0=\startpage

\gt\hfill      
\beginpicture
\setcoordinatesystem units <0.33truein, 0.33truein> point at 2.2 0.9
\setplotsymbol ({$\cal G$})
\plotsymbolspacing=9truept
\circulararc 315 degrees from 0 1 center at 0 0
\setplotsymbol ({$\cal T$})
\circulararc 315 degrees from 1 -1 center at 1 0
\endpicture
%
\break
{\small\ifx\thesamplenumber\relax 
Volume \else Sample
\fi\thevolumenumber\ (\thevolumeyear)
\startpage--\finishpage\nl
Published: \publishdate}
\vglue 0.5truein plus 0.4fil minus 0.1truein

{\parskip=0pt\leftskip 0pt plus 1fil\def\\{\par\smallskip}{\ifplaintex\large
\else\Large\fi\bf\thetitle}\par\medskip}   

\vglue 0pt plus 0.1fil 

{\parskip=0pt\leftskip 0pt plus 1fil\def\\{\par}{\sc\theauthors}
\par\medskip}

\vglue 0pt plus 0.1fil 

{\small\parskip=0pt\let\newline\\
{\leftskip 0pt plus 1fil\def\\{\par}{\sl\theaddress}\par}
\expandafter\ifx\theemail\relax    
\relax\else\vglue 5pt plus 0.02fil minus 2pt\def\\{\stdspace{\rm 
and}\stdspace} 
\cl{Email:\stdspace\tt\theemail}\fi
\ifx\theurl\relax                  
\relax\else\vglue 5pt plus 0.02fil minus 2pt\def\\{\stdspace{\rm 
and}\stdspace}
\cl{URL:\stdspace\tt\theurl}\fi\par}

\vglue 7pt plus 0.3fil minus 3pt

{\bf Abstract}
\vglue 5pt plus 0.1fil minus 2pt

\theabstract

\vglue 7pt plus 0.3fil minus 3pt

{\bf AMS Classification numbers}\quad Primary:\quad \theprimaryclass

Secondary:\quad \thesecondaryclass

\vglue 5pt plus 0.3fil minus 2pt

{\bf Keywords}\quad \thekeywords

\vglue 10pt plus 0.5fil minus 5pt

{\small  Proposed: \theproposer\hfill Received: \receiveddate\nl
Seconded: \theseconders\hfill 
\ifx\reviseddate\relax                         
Accepted: \accepteddate                        
\else
Revised: \reviseddate                          
\fi}
\eject
}       

\let\maketitlepage\maketitlep
\let\maketitle\maketitlepage


\font\phead=cmsl9 scaled 950
\font\lhead=cmsl9 scaled 1050
\font\pnum=cmbx10 scaled 913
\font\lnum=cmbx10 
\font\pfoot=cmsl9 scaled 950
\font\lfoot=cmsl9 scaled 1050
\ifplaintex
\headline{\vbox to 0pt{\vskip -4.5mm\line{\small\phead\ifnum
\count0=\startpage ISSN 1364-0380 (on line)
1465-3060 (printed) \hfill {\pnum\folio}\else\ifodd\count0\def\\{ }%
\ifx\theshorttitle\relax\thetitle\else\theshorttitle\fi\hfill{\pnum\folio}
\else\def\\{ and }{\pnum\folio}\hfill\ifx\theshortauthors\relax\theauthors
\else\theshortauthors\fi\fi\fi}\vss}}
\footline{\vbox to 0pt{\vglue 0mm\line{\small\pfoot\ifnum\count0=\startpage
\copyright\ \gtp\hfill\else
\gt, Volume \thevolumenumber\ (\thevolumeyear)\hfill\fi}\vss
}}
\else
\makeatletter
\def\@oddhead{{\small\lhead\ifnum\count0=\startpage ISSN 1364-0380 (on line)
1465-3060 (printed) \hfill {\lnum\number\count0}\else\ifodd\count0
\def\\{ }\ifx\theshorttitle\relax \thetitle \else\theshorttitle\fi\hfill
{\lnum\number\count0}\else\def\\{ and }{\lnum\number\count0}
\hfill\ifx\theshortauthors\relax 
\theauthors\else\theshortauthors\fi\fi\fi}}\def\@evenhead{\@oddhead}
\def\@oddfoot{\small\lfoot\ifnum\count0=\startpage\copyright\ \gtp\hfill\else
\gt, Volume \thevolumenumber\ (\thevolumeyear)\hfill\fi}
\def\@evenfoot{\@oddfoot}
\makeatother
\fi


\newwrite\gtoutfile
\long\gdef\makeheadfile{  
{\def\\{, }\def\s{ }
\immediate\openout\gtoutfile head.xxx
\immediate\write\gtoutfile{To: math@arxiv.org}
\immediate\write\gtoutfile{Subject: put or rep NNNNN:pppp}
\immediate\write\gtoutfile{--text follows this line--}
\immediate\write\gtoutfile{Proxy-for: \ifx\theasciiauthors\relax
\theauthors\else\theasciiauthors\fi\s<\ifx\theasciiemail\relax\theemail\else\theasciiemail\fi>}
\immediate\write\gtoutfile{\noexpand\\}
\immediate\write\gtoutfile{Authors: \ifx\theasciiauthors\relax
\theauthors\else\theasciiauthors\fi}
{\def\\{ }\immediate\write\gtoutfile{Title: \ifx\theasciititle\relax
\thetitle\else\theasciititle\fi}}
\immediate\write\gtoutfile{Subj-class: GT or SG or MG etc}
\immediate\write\gtoutfile{MSC-class: \theprimaryclass\ifx\thesecondaryclass\relax\else, \thesecondaryclass\fi}
\immediate\write\gtoutfile{Journal-ref: Geom. Topol. \thevolumenumber
(\thevolumeyear) \startpage-\finishpage}
\immediate\write\gtoutfile{Comments: Published by Geometry and Topology at}
\immediate\write\gtoutfile{\s\s http://www.maths.warwick.ac.uk/gt/GTVol\thevolumenumber/paper\thepapernumber.abs.html}
\immediate\write\gtoutfile{\noexpand\\}
\immediate\write\gtoutfile{}
\ifx\theasciiabstract\relax
\immediate\write\gtoutfile{\theabstract}\else
\immediate\write\gtoutfile{\theasciiabstract}\fi
\immediate\write\gtoutfile{}
\immediate\write\gtoutfile{\noexpand\\}
\immediate\write\gtoutfile{}
\immediate\closeout\gtoutfile}}  

\def\maketitlepage{\maketitlep\makeheadfile}
\let\maketitle\maketitlepage

\lognumber{248}

\volumenumber{6}
\papernumber{17}
\volumeyear{2002}
\pagenumbers{495}{521}
\received{20 April 2002}
\revised{19 November 2002}
\accepted{19 November 2002}
\published{22 November 2002}
\proposed{David Gabai}
\seconded{Jean-Pierre Otal, Joan Birman}

\input epsf

\let\bold\Bbb
\def\tag#1${\eqno{(#1)}$}

\def\frac#1#2{{#1\over#2}}

\def\S{section }

\reflist

\key{Be} {\bf A Beardon}, {\it The geometry of discrete groups}, 
Springer--Verlag, Berlin--New York (1983)

\key{Bo1} {\bf F Bonahon}, {\it Bouts des vari\'et\'es hyperboliques 
de dimension $3$},  Ann. of Math. 124 (1986) 71--158

\key{Bo2} {\bf F Bonahon}, 
{\it Earthquakes on Riemann surfaces and on measured geodesic laminations}, 
Trans. Amer. Math. Soc. 330 (1992) 69--95 

\key{Bu} {\bf P Buser}, {\it Geometry and spectra of compact Riemann surfaces},
Birkh\"auser, Boston (1992)

\key{De} {\bf M Dehn}, {\it Papers on group theory and topology}, J. Stillwell 
(editor),  Springer--Verlag, Berlin--New York (1987)

\key{FLP} {\bf A Fathi}, {\bf F Laudenbach}, {\bf V Poenaru}, {\it 
Travaux de Thurston sur les surfaces}, Ast\'erisque 66--67, Soci\'et\'e
Math\'ematique de France (1979)

\key{HM}  {\bf J Hubbard}, {\bf H Masur}, {\it Quadratic differentials and 
foliations}, Acta Math. 142 (1979) 221--274

\key{IT} {\bf Y Imayoshi}, Y. and {\bf M Taniguchi}, {\it An introduction 
to Teichm\"uller spaces}, Translated and revised from the Japanese by
the authors, Springer--Verlag, Tokyo (1992)

\key{Lu} {\bf F Luo}, {\it Simple loops on surfaces and their 
intersection numbers}, preprint (1997)

\key{LS} {\bf F Luo}, {\bf R Stong}, {\it Dehn--Thurston coordinates 
of curves on surfaces}, preprint (2002)

\key{Pa} {\bf A Papadopoulos},
{\it On Thurston's boundary of Teichm\"uller space and the extension
of earthquakes} Topology Appl. 41 (1991) 147--177

\key{PH} {\bf R Penner}, {\bf J Harer},
{\it Combinatorics of train tracks}, Annals of Mathematics Studies,
125, Princeton University Press, Princeton, NJ (1992)

\key{Th1}  {\bf W Thurston}, {\it On the geometry and  dynamics of 
diffeomorphisms of surfaces}, Bull. Amer. Math. Soc. 19 (1988) 417--438

\key{Th2} {\bf W Thurston}, {\it Geometry and topology of 3--manifolds}, 
Princeton University lecture notes (1976)

\endreflist

\title{Lengths of simple loops on surfaces\\with hyperbolic 
metrics}  
\author{Feng Luo\\Richard Stong}    
 
\address{Department of Mathematics, Rutgers University\\New Brunswick, 
NJ 08854, USA}

\secondaddress{Department of Mathematics, Rice University\\Houston, TX 77005,
USA}

\asciiaddress{Department of Mathematics, Rutgers University\\New Brunswick, 
NJ 08854, USA\\and\\Department of Mathematics, Rice University\\Houston, 
TX 77005,
USA}

\email{fluo@math.rutgers.edu}

\secondemail{stong@math.rice.edu}

\asciiemail{fluo@math.rutgers.edu, stong@math.rice.edu}

\abstract
Given a compact orientable surface of negative Euler characteristic,
there exists a natural pairing between the Teichm\"uller space of the
surface and the set of homotopy classes of simple loops and arcs. The
length pairing sends a hyperbolic metric and a homotopy class of a
simple loop or arc to the length of geodesic in its homotopy class. We
study this pairing function using the Fenchel--Nielsen coordinates on
Teichm\"uller space and the Dehn--Thurston coordinates on the space of
homotopy classes of curve systems.  Our main result establishes
Lipschitz type estimates for the length pairing expressed in terms of
these coordinates.  As a consequence, we reestablish a result of
Thurston--Bonahon that the length pairing extends to a continuous map from the
product of the Teichm\"uller space and the space of measured
laminations.
\endabstract                  

\asciiabstract{Given a compact orientable surface of negative Euler 
characteristic, there exists a natural pairing between the
Teichmueuller space of the surface and the set of homotopy classes of
simple loops and arcs. The length pairing sends a hyperbolic metric
and a homotopy class of a simple loop or arc to the length of geodesic
in its homotopy class. We study this pairing function using the
Fenchel-Nielsen coordinates on Teichmueller space and the
Dehn-Thurston coordinates on the space of homotopy classes of curve
systems.  Our main result establishes Lipschitz type estimates for the
length pairing expressed in terms of these coordinates.  As a
consequence, we reestablish a result of Thurston-Bonahon that the
length pairing extends to a continuous map from the product of the
Teichmueller space and the space of measured laminations.}

\keywords{Surface, simple loop, hyperbolic metric, Teichm\"uller space}
\asciikeywords{Surface, simple loop, hyperbolic metric, Teichmueller space}

\primaryclass{30F60}
\secondaryclass{57M50, 57N16}

\maketitlepage

\section{Introduction}

\rk{1.1}Given a compact orientable surface  of negative Euler 
characteristic, there exists a natural length pairing between the
Teichm\"uller space of the surface and the set of homotopy classes of
simple loops and arcs. The length pairing sends a hyperbolic metric
and a homotopy class of a simple loop or arc to the length of the geodesic
in its homotopy class. In this paper, we study this pairing function
using the Fenchel--Nielsen coordinates on Teichm\"uller space and the
Dehn--Thurston coordinates on the space of homotopy classes of curve
systems. Our main result, theorem 1.1, establishes Lipschitz type
estimates for the length pairing expressed in terms of these
coordinates. As a consequence, we give a new proof of a result of
Thurston--Bonahon ([\ref{Th1}], see [\ref{Bo1}, proposition 4.5] for a proof) that
the length pairing extends to a continuous map from the product of the
Teichm\"uller space and the space of measured laminations to the real
numbers so that the extension is homogeneous in the second coordinate.

\rk{1.2} Let $F$ be a compact connected orientable surface with possibly 
non-empty
boundary and negative Euler characteristic. By a hyperbolic metric
on the surface $F$ we mean a Riemannian metric of curvature $-1$ on the 
surface $F$ so that its boundary components are geodesics.
The Teichm\"uller space $T(F)$ is the space of all isotopy classes of
hyperbolic metrics on the surface. Recall that two hyperbolic metrics 
are
\it isotopic \rm if there is an isometry between the two metrics which 
is
isotopic to the identity. Following M. Dehn [\ref{De}], a \it curve system 
\rm in
the surface $F$ is a compact proper 1--dimensional submanifold so that 
each
of its circle components is not null homotopic and not homotopic into 
the
boundary $\partial F$ of $F$ and each of its arc component is not 
homotopic
into $\partial F$ relative to its endpoints. We denote the set of all
homotopy classes (or equivalently isotopy classes) of curve systems on
$F$ by $CS(F)$ and
call it \it the space of curve systems. \rm By a basic fact from 
hyperbolic
geometry, for any hyperbolic metric $d$ on $F$ and
any homotopically non-trivial simple loop or arc $s$ in $F$, there is a
unique shortest $d$--geodesic $s^*$ homotopic (and isotopic) to $s$. One
defines the length of the homotopy class $[s]$, denoted by $l_d([s])$ 
(or
$l_{[d]}([s])$ since it depends only on the class $[d]\in T(F)$), to be 
the
$d$--length of the geodesic $s^*$. This length pairing
extends naturally to a map $T(F) \times CS(F) \to \bold R$, still 
denoted
by $l_d([s])$. 

Our goal is to understand this length pairing using parametrizations of
$T(F)$ and $CS(F)$. To this end, let us recall the Fenchel--Nielsen 
coordinates
on Teichm\"uller space and Dehn--Thurston coordinates on the
space of curve systems. The definition of these two coordinates depends 
on
the choice of a hexagonal decomposition on the surface (see \S2.2).
Fix such a decomposition on a surface of genus $g$ with $r$ boundary
components, we obtain a parametrization (the Fenchel--Nielsen 
coordinates)
of the Teichm\"uller space $FN\co  T(F) \to  R$ where
$R =(\bold R_{>0} \times \bold R)^{3g-r+3} \times \bold R_{>0}^r$ and 
a
parametrization (the  Dehn--Thurston coordinates) $DT\co  CS(F) \to Z$ 
where
$Z = ((\bold Z \times \bold Z)/\pm)^{3g-r+3} \times \bold Z_{\geq 
0}^r$.
(See \S2 and \S3 for details). Here $\bold R_{>0}$ and $\bold Z_{>0}$ denote the sets of
positive real numbers and positive integers respectively. Note that $FN$ is a homeomorphism
and $DT$ is an (homogeneous) injective map.

We introduce a metric on the space $Z$ as follows. 
The metric on $(\bold Z \times \bold Z)/\pm$ is defined to be
$|(x_1, y_1)-(x_2, y_2)|=\min\{|x_1+x_2|+|y_1 + y_2|, |x_1-x_2|+|y_1 - 
y_2|\}$.
The metric on $\bold Z_{>0}$ is the standard metric and the metric on 
$Z$ is
the product metric. The length $|x|$ of
$x = ([x_1, t_1],..., [x_N, t_N], x_{N+1}, ..., x_{N+r}) \in Z$ is
$\sum_{i=1}^{N+r} |x_i| + \sum_{j=1}^N |t_j|$ where $N=3g+r-3$.

For $x=(x_1, t_1,..., x_N,t_N,x_{N+1},\dots , x_{N+r})$ and\nl
\hbox{}\hglue 1.8in $y = (y_1, s_1, ..., y_N, s_N, y_{N+1},\dots , y_{N+r})$ in $R$,
let 
$$\eqalign{D(x,y) = \smash{\sum_{i=1}^N} \min\{x_i, y_i\} |t_i -s_i|& +\cr
 (5 \max_i \{|t_i|,|s_i|\}+7)&
\sum_{i=1}^{N+r} |\log \sinh(x_i/2) - \log \sinh ( y_i/2)|.}$$
Note that this $D\co  R \times R \to \bold R$ is continuous and satisfies
$D(x,y)>0$ if $x \neq y$, but it is not a metric on $R$. Define
$$|x| = \sum_{i=1}^N ( x_i  + 1/x_i + x_i |t_i|)   + \sum_{j=N+1}^{N+r}
(x_j + 1/x_j) + (N+r)\log 2.$$
Here $x_i$ is the length of the $i$-th decomposing
loop in the metric and $x_i t_i$ is the twisting length. The number
$2 \pi t_i$ measures the angle of twisting at the $i$-th decomposing
loop.

Our main theorem is the following.

\proclaim{Theorem 1.1} Suppose $F$ is a compact orientable surface with
possibly non-empty boundary components and the surface $F$ has a fixed
hexagonal decomposition. Let
$FN\co T(F)\to \bold (\bold R_{>0} \times \bold R)^{3g-r+2} \times \bold 
R_{>0}^r$
and
$DT\co CS(F)\to ((\bold Z \times \bold Z)/\pm)^{3g-r+2} \times \bold 
Z_{\geq 0}^r$
be the Fenchel--Nielsen coordinate and the Dehn--Thurston coordinate 
associated
to the hexagonal decomposition. Then for any $[a], [b]$ in $CS(F)$ and 
any two
hyperbolic metrics $[d_1], [d_2]$ in $T(F)$, the following inequalities 
hold.
$$|l_{d_1}([a]) - l_{d_1}([b])| \leq 3|FN(d_1)||DT([a])-DT([b])|, \tag 
1.1$$
and
$$|l_{d_1}([a])-l_{d_2}([a])| \leq  4D( FN(d_1), FN(d_2)) |DT([a])|. 
\tag 1.2$$\endproc

As a consequence, we give a new proof of the following result of Thurston--Bonahon (see [\ref{Bo1}] for the first
published proof).

\proclaim{Corollary 1.2}{\rm ([\ref{Th1}], [\ref{Bo1}])}\qua 
The hyperbolic length function extends to a continuous map from $T(F)
\times ML(F) \to \bold R$ where $ML(F)$ is the space of measured
laminations on the surface $F$. Furthermore, the extension also
satisfies the inequalities (1.1) and (1.2). \endproc

\rk{1.3}One of the main ingredients used in the proof is the following 
elementary
geometric fact about right-angled hyperbolic hexagons (see theorem 5.2 
in \S5).
Let $H_{x,a,b}$ be a right-angled hyperbolic hexagon whose side lengths are 
(reading
from counterclockwise): $a,z,x,y,b,w$. Let $S_{\lambda, \mu}$ be the 
length of
a geodesic segment in $H_{x,a,b}$ joining any two sides of the hexagon so that
the endpoints of the segment cut the sides into two intervals of 
lengths
$\lambda t, (1-\lambda)t$ and $\mu s$, $(1-\mu) s$. Then if we fix 
$a,b$,
$\lambda$, $\mu$ and  let $x$ vary, 
the length $S_{\lambda,\mu}$ satisfies
$$ \left|\frac{ d S_{\lambda, \mu}}{dx} \right| \leq 4 \coth x.$$
In particular, this implies that,
$$ |S_{\lambda, \mu}(x) - S_{\lambda, \mu}(x')| \leq
4 | \log \sinh(x) - \log \sinh(x')|.$$

\midinsert\centerline{\epsfxsize 1.5in \epsfbox{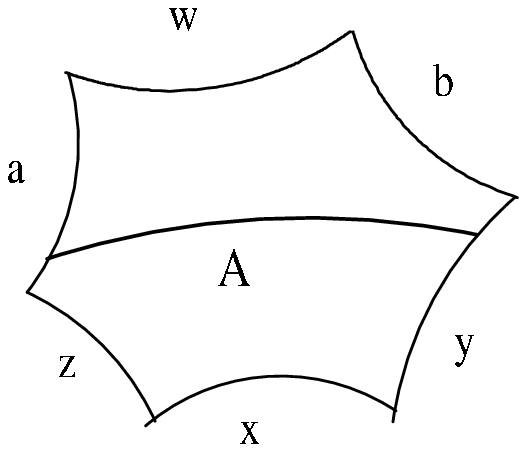}}
\vglue 6pt
\centerline{\small Figure 1.1}
\endinsert

\rk{1.4} The paper is organized as follows. In section 2, we recall some of 
the
known facts about the curve systems and the results obtained in [\ref{LS}]. 
In
particular, we will recall the notion of the hexagonal decompositions 
of the
surface and the Dehn--Thurston coordinates on the space of curve 
systems. In
section 3, we will recall the Fenchel--Nielsen coordinates of hyperbolic
metrics. The main theorem 1.1 will be proved in section 4. In section 
5, we
establish two simple facts on hyperbolic right angled hexagon used in 
the
proof.

The work is supported in part by the NSF. 

\section{Dehn--Thurston coordinates of curve systems}

We will recall the Dehn--Thurston coordinates on $CS(F)$ in this 
section. The
basic ingredient to set up the coordinate is the colored hexagonal
decomposition of a surface which is defined in subsection 2.1 below. 
Unless mentioned otherwise, we will assume in this section that the 
surface
$F$ is oriented with negative Euler characteristic.

\sh{2.1\qua Notation and conventions}

We shall use the following notations and conventions. Let $F = F_{g, r}$ be 
the
orientable compact surface of genus $g$ with $r \geq 0$ boundary 
components.
The interior of a surface $F$ will be denoted by $int(F)$. All 
subsurfaces
in an oriented surface have the induced orientation. We will always 
draw
oriented surface so that its orientation is the right-hand orientation 
on the front face of the surface that we see. 

A \it curve system \rm on $F$ is a proper 1--dimensional submanifold $s$ 
in
$F$ so that no circle component of $s$ is  null homotopic or 
homotopic into the boundary of the surface $F$ and no arc component 
of
$s$ is  null homotopic relative to the boundary. If $s$ is a proper
submanifold of a surface, we use $N(s)$ to denote a small tubular 
neighborhood
of $s$. The isotopy class of a submanifold $s$ is denoted by $[s]$. If 
$a$ 
and $b$ are isotopic submanifolds we will write $a \cong b$. If $a, b$ 
are
two proper 1--dimensional submanifolds, we will use $I(a, b)$, $I([a], 
b)$ or
$I(a, [b])$ to denote the geometric intersection number
$I([a], [b]) = \min\{ |a' \cap b'|: a \cong a', b \cong b'\}$. 
Here $|X|$ denoted the cardinal of a set $X$.
When a curve system $a$ is written as a union $a_1 \cup ...\cup a_n$, 
it is
understood that each $a_i$ is a union of components of $a$. Let $2 
\bold Z$
be the set of even integers. All hyperbolic metrics on compact surfaces
are assumed to have geodesic boundary. Also if $d$ is a hyperbolic 
metric and
$a$ is a curve system, we use $l_d(a)$ to denote the length of $a$ in 
the
metric $d$. The length of the isotopy class $[a]$ is defined to be
$\inf \{ l_d(a')| a' \cong a\}$ and is denoted by $l_d([a])$.

Fix an orientation on the surface $F$. 
Let us recall the concept of multiplication of two curve systems in
$CS(F)$ (see [\ref{Bo2}], [\ref{Pa}] and [\ref{Lu}], the notation was first 
introduced in
[\ref{Pa}], [\ref{Bo2}] as the earthquakes in the space of measured laminations).
Given $\alpha$ and $\beta$ in $CS(F)$, take $a \in \alpha$ and
$b \in \beta$ so that $|a \cap b|=$ $I(\alpha, \beta)$. If $\alpha$ and
$\beta$ are disjoint, we define $\alpha \beta$ to be $[a \cup b]$. If
$I(\alpha, \beta)$ $> 0$, then $\alpha \beta$ is defined to be the 
isotopy
class of the 1--dimensional submanifold $ab$ obtained by resolving all
intersection points in $a \cap b$ from $a$ to $b$. Here by the 
resolution
from $a$ to $b$ we mean the following surgery. At each point $p \in a 
\cup b$,
fix any orientation on $a$. Then use the orientation of the surface to
determine an orientation of $b$ at $p$. Finally resolve the singularity 
at $p$
according to the orientations on $a$ and $b$. One checks easily that 
this is
independent of the choice of orientation on $a$. See figure 2.1. If $a$ 
is a
curve system and $k$ is a positive integer, then the collection of $k$
parallel copies of $a$ is denoted by $a^k$. We use $[a]^k$ to denote 
$[a^k]$.
If $k$ is a negative integer, we denote $[a]^k[b]$ by $[b] [a]^{-k}$ 
and
$[b][a]^k$ by $[a]^{-k}[b]$.


\midinsert
\centerline{\epsfxsize 3.5in\epsfbox{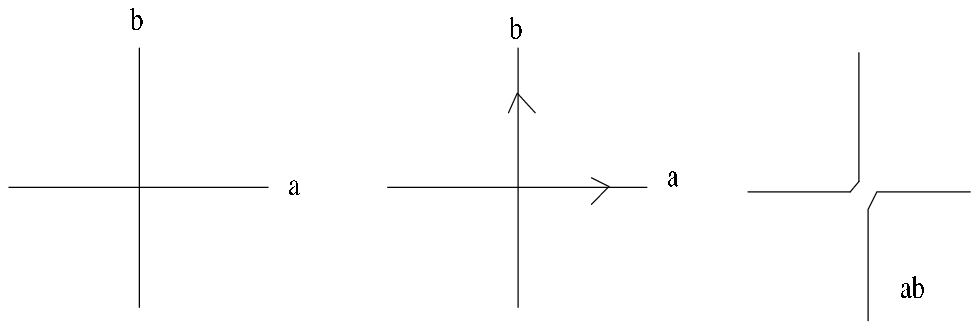}}
\vglue 6pt
\centerline{\small Figure 2.1}  
\endinsert

The following useful property follows from the definition.

\proclaim{Lemma 2.1} {\rm(Triangle inequality)}\qua Suppose $a$ is a curve 
system
without arc components and $b$ is a curve system. Fix a hyperbolic metric $d$ 
on
the surface $F$. Then the hyperbolic lengths satisfy
$$ | l_d( [ab]) - l_d([b]) | \leq l_d([a]),$$
$$ |l_d([ba]) - l_d([b]) | \leq l_d([a]).\leqno{{\sl and}}$$ \endproc

Indeed, by the definition of resolutions and taking all components of 
$a$ and
$b$ to be geodesics, one sees that $l_d([ab]) \leq l_d([a]) + l_d([b])$
(this inequality also holds for curve systems $a$ with arc components). 
To see
the inequality $l_d([b]) \leq l_d([ab]) + l_d([a])$, we use the 
cancelation
property of the multiplication ([\ref{LS}] theorem 2.4(4)) that
$(ab)a  \cong b \cup c^2$ where $c$ consists of those components of $a$ 
which
are disjoint from $b$. Thus
$l_d([b])  \leq l_d([b \cup c^2]) = l_d( [(ab)a]) \leq l_d([ab]) + 
l_d([a])$.
This proves the lemma.

A curve system $s$ on $F$ is called a \it 3--holed sphere decomposition 
\rm if 
(1) each component of $s$ is a circle and (2) all components of $F -s$ 
are
3--holed spheres.  This implies that $s$ contains $3g+r-3$ many 
components when
$F=F_{g,r}$.

By a \it hexagonal decomposition \rm of the 3--holed sphere $F_{0,3}$, 
we mean
a curve system $b$  on $F_{0,3}$ so that $b$ contains exactly three arc
components joining different boundary components in $F_{0,3}$. See 
figure
2.2(a). We call each component of $F_{0,3}-b$ a \it hexagon. \rm \it A
colored hexagonal decomposition \rm of an orientable compact surface 
$F$ is a
triple $(p, b, col)$ where $p, b$ are curve systems and $col$ is a 
coloring
so that (1) $p$ is a 3--holed sphere decomposition, (2) for each 
component
$F'$ of $F-p$, the intersection $b \cap F'$ is a hexagonal 
decomposition of
the 3--holed sphere, (3) one can color the components of $F -p\cup b$ 
into
red and white so that there is exactly one red hexagon in each 
component
of $F-p$ and the red hexagons join only red hexagons crossing $p$. The
triple $(p, b, col)$ is also called a \it marking \rm on the surface 
$F$.

\rk{2.2}The classification of the curve systems on the 3--holed sphere 
$F_{0,3}$
is well known. Suppose the  boundary components of the 3--holed sphere
$F_{0,3}$ are $\partial_1, \partial_2, \partial_3$. Then each
$[a] \in CS(F_{0,3})$ is determined uniquely by $DT([a]) =(x_1, x_2, 
x_3)$
where $x_i =I(a, \partial_i)$. Furthermore the map
$DT\co CS(F_{0,3}) \to$\break $\{(x_1, x_2, x_3) \in \bold Z_{\geq 0}^3 | x_1 + 
x_2 + x_3  \in 2 \bold Z\}$
is a bijection. These are the Dehn--Thurston coordinates for the 3--holed 
sphere.
The curve systems with coordinates $(x_1, x_2, x_3)$ are shown in 
figure 2.2(b).

\midinsert
\centerline{\epsfxsize 3.5in \epsfbox{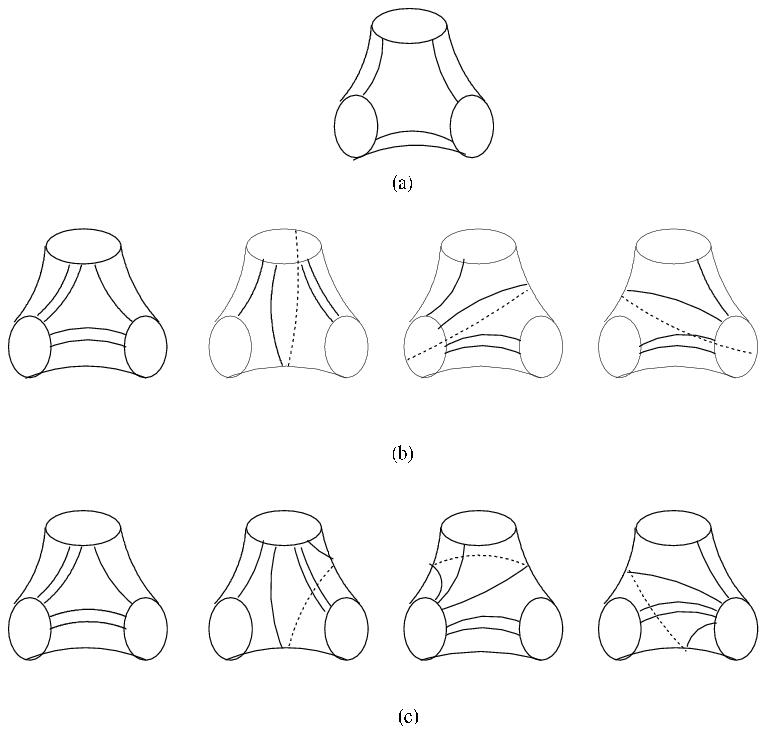}}
\vglue 6pt 
\centerline{\small Figure 2.2} 
\endinsert

If we fix a colored hexagonal decomposition $b=b_1\cup b_2 \cup b_3$ of 
the
oriented surface $F_{0,3}$, then each $[a] \in CS(F_{0,3})$ has a \it
standard representative \rm with respect to the hexagonal 
decomposition. It
is defined as follows. We assume that $b_i$ is disjoint from 
$\partial_i$.
Take a curve system $a$ in $F_{0,3}$. Its standard representative is a 
curve
system $a' \cong a$ so that each component of $a'$ is standard. Here an 
arc
$s$ is \it standard \rm if either it lies entirely in the red-hexagon 
or if
$\partial s \subset \partial_i$, then $\partial s$ is in the 
red-hexagon and
$| s \cap (b_1 \cup b_2 \cup b_3)| = 2 = |s \cap (b_i \cup b_j)|$ so 
that
the cyclic order of the sets $(s\cap \partial_i, s \cap b_i, s \cap 
b_j)$
in the boundary of the red-hexagon coincides with the induced 
orientation
from the red-hexagon. For instance the standard representatives of
the curve systems with coordinates $(x_1, x_2, x_3)$ are shown in 
figure
2.2(c) where the red-hexagon is the front hexagon in figure 2.2(a).

Fix a marking  $(p_1 \cup \dots \cup p_{3g+r-3}, b, col)$ on an 
oriented
surface $F = F_{g,r}$. The Dehn--Thurston coordinates of $[a]$ in 
$CS(F)$ is
a vector in $(\bold Z^2/\pm)^{3g+r-3} \times \bold Z^r_{\geq 0}$ 
defined as
follows. Express the class $[a]$ as
$$ [a] = [p_1^{t_1} \dots p_{3g+r-3} ^{t_{3g+r-3}}] [a_{zt}]$$
where $t_i \in \bold Z$ so that if $I(a, p_i) =0$ then $t_i \geq 0$ 
and $a_{zt}$ is a curve system so that its restriction to each 3--holed 
sphere
component of $F -p_1 \cup \dots\cup p_{3g+r-3}$ is a standard curve 
system
with respect to the red hexagon. Then the Dehn--Thurston coordinate of 
$[a]$ is
$$DT([a])=([x_1,t_1],\dots,[x_{3g+r-3},t_{3g+r-3}],x_{3g+r-2},\dots,x_{3g+2r-3})$$
where $x_i = I(a, p_i)$ and $p_{3g+r-3+j} = \partial_j F$. Note that
$I(a, p_i) = I(a_{zt}, p_i)$ and the twisting coordinates $t_i(a_{zt})$ of $a_{zt}$ 
are zero.
We sometimes use $x_i(a)$ and $t_j(a)$ to denote the coordinates $x_i$ 
and
$t_j$ of the curve systems $a$. It is shown in [\ref{LS}] (proposition 2.5) 
that
this is well defined. For $[s] \in CS(F)$ and $k \bold Z_{>0}$, let
$ [s]^k = [s^k]$ be the isotopy class of $k$--parallel copies of $s$.

\proclaim{Proposition 2.2}
 The Dehn--Thurston coordinate is a bijection 
$$\eqalign{
DT\co CS(F)&\to\cr 
\{ ([x_1,t_1],&
\dots,[x_{3g+r-3},t_{3g+r-3}],x_{3g+r-2},\dots,
x_{3g+2r-3})\in (\bold Z^2/\pm)^{2g+r-3}\cr \times (\bold Z_{\geq 0})^r \mid
\hbox{ if } &p_i, p_j \hbox{ and }p_k\hbox{ bound a 3--holed sphere, then }
x_i + x_j + x_k \in 2 \bold Z\}.}$$
Furthermore, $DT([a]^k) = k DT([a])$ 
for
$k \in \bold Z_{\geq 0}$.\endproc

\sh{2.3\qua The main idea of the proof of theorem 1}

We sketch the proof of the inequality (1.1) in the main theorem 1.1 in
this subsection. First of all, by homogeneity $l_d([a^2]) = 2l_d([a])$
and $DT(a^2)= 2DT(a)$, hence it suffices to prove (1.1) for classes
$[a], [b]$ so that $DT(a) =u$ and $DT(b) = v$ are \it even vectors\rm,
ie, all $x_i$ and $t_j$ coordinates of them are even integers. Now
given any two even vectors $u$ and $v$ in $Z$ with distance $|u-v| =
2n$ there exists a sequence of $n+1$ even vectors $u_0 = u, u_1,
\dots, u_n = v$ so that $|u_i - u_{i+1}| = 2$. On the other hand, by
proposition 2.2, each even vector $u_i$ is the image $DT(a_i)$ for
some $[a_i] \in CS(F)$. Thus by interpolation, it suffices to prove
inequality (1.1) for classes $[a]$ and $[b]$ so that $DT(a)$ and
$DT(b)$ are even vectors of distance two apart. This means that the
Dehn--Thurston coordinates of $[a]$ and $[b]$ are the same except at
one $x_i$-- or $t_j$--coordinate where they differ by $2$.  If one of
their twisting coordinates differs by 2, say $t_i(a) = t_i(b)+2$, then
$[a] = [p_i^2 b]$ by definition. Thus, by the triangle inequality
(Lemma 2.1), we have $|l_d([a])-l_d([b])|\leq
l_d([p_i^2])=2l_d([p_i])\leq |FN(d)||DT(a)-DT(b)|$.  If their
intersection number coordinates differ by two, say $x_i(a) = x_i(b) +
2$, for some $i$ with $1 \leq i \leq 3g+r -3$, then we prove in
[\ref{LS}] (proposition 4.3) that $[a] = \delta_1 ... \delta_s [b]
\delta_{s+1} ... \delta_t$ where $t
\leq 5$
and the $\delta_i$'s are quite simple. In fact, we show that these 
simple
loops $\delta_i$'s satisfy 
$$ \sum_{i=1}^t l_d(\delta_i) \leq 6 |FN(d)|.$$
Thus by the triangle inequality (lemma 2.1),
$|l_d([a])-l_d([b])|\leq \sum_{i=1}^t l_d(\delta_i)\leq 
6|FN(d)|=3|FN(d)||DT(a)-DT(b)|.$
If their intersection number coordinates differ by two $x_i(a) = x_i(b) 
+2$ for
some $i$ with $i \geq 3g+r-2$, then doubling the surface across its 
boundary
reduces to the previous case.

This shows that the main issue is to understand the effect of changing 
some
intersection coordinate $x_i$ by $2$. This will be addressed in the 
following
subsections.

\rk{2.4}
We will recall the results obtained in [\ref{LS}] concerning the
change of $x_i$ coordinates by 2. Suppose $(p_1 \cup \dots \cup
p_{3g+r-3}, b, col)$ is a marking on an oriented surface $F$, and $DT$
is the associated Dehn--Thurston coordinate. Let $[a]$ and $[b]$ be two
isotopy classes of curve systems so that their twisting coordinates
$t_j(a)$ and $t_j(b)$ are the same and their intersection coordinates
agree except for the $i$-th which satisfies $x_i(a) = x_i(b) + 2$. We
will find a surgery procedure converting $a$ to $b$.  There are three
cases to be discussed. In the first case, the corresponding
decomposing simple loop $p_i$ is adjacent to only one 3--holed sphere
component of $F-p$ and $p_i$ is not in $\partial F$. In the second
case, the simple loop $p_i$ is adjacent to two different components of
$F-p$. In the last case, $p_i$ is a boundary component of the surface
$F$.

\midinsert
\centerline{\epsfbox{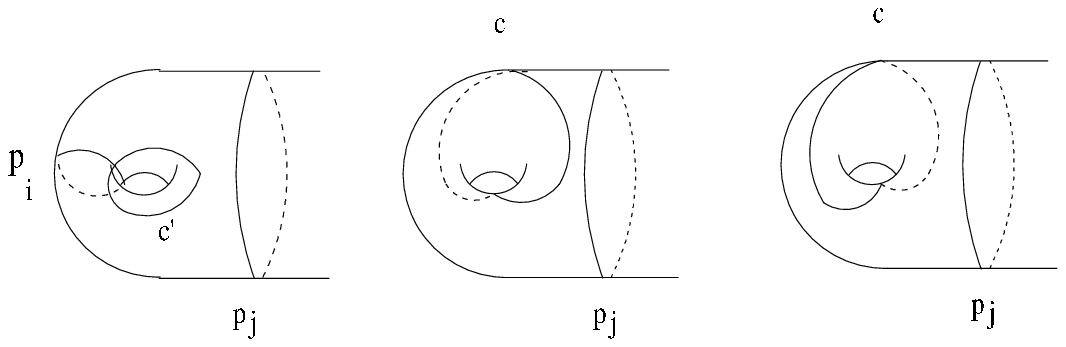}}
\vglue 6pt
\centerline{\small 
Figure 2.3: Here $c'$ is the simple loop with zero twisting coordinate.}
\centerline{\small 
The loop $c$ is obtained from $c'$ by a Dehn twist along $p_i$.}
\endinsert

The following two results were obtained in [\ref{LS}] (propositions 4.2 and 
4.3).

\proclaim{Proposition 2.3} {\rm([\ref{LS}], proposition 4.2)}\qua
In the first case that $p_i$ is adjacent to only one 3--holed sphere,
suppose $p_j$ is the simple loop bounding the 1--holed torus which
contains $p_i$. Then
$$ a \cong p_j^{e_1} c^{e_2} b$$
where $e_1, e_2 \in \{0, \pm 1, \pm 2\}$ and $c$ is one of the two 
simple loops
with Dehn--Thurston coordinates
$([0,0],\dots, [0,0], [1,\! \pm 1], [0,0],\dots, [0,0],0,\dots,0)$
(the non\-zero coordinates are $x_i$ and $t_i$). See figure 2.3. \endproc

\proclaim{Proposition 2.4} {\rm([\ref{LS}], proposition 4.3)}\qua
In the second case that $p_i$ is adjacent to two 3--holed spheres,
suppose $p_{i_1},\dots, p_{i_4}$ are the simple loops bounding the
4--holed sphere containing $p_i$ and $p_i, p_{i_1}, p_{i_2}$ bound a
3--holed sphere.  Then
$$a \cong p_{i_1}^{s_1}\dots p_{i_4}^{s_4} c^e b$$
where $e \in \{\pm 1\}$, $|s_1| + |s_2| \leq 2$, $|s_3|+ |s_4| \leq 2$ 
and
$c$ is a simple loop in the 4--holed sphere whose Dehn--Thurston 
coordinates are
$DT(c)=([0,0],\dots ,[2, t],$
$\dots ,[0,0],0,\dots ,0)$ so that $|t|\leq 
2$. See figure 2.4. \endproc

\midinsert
\centerline{\epsfxsize.9\hsize\epsfbox{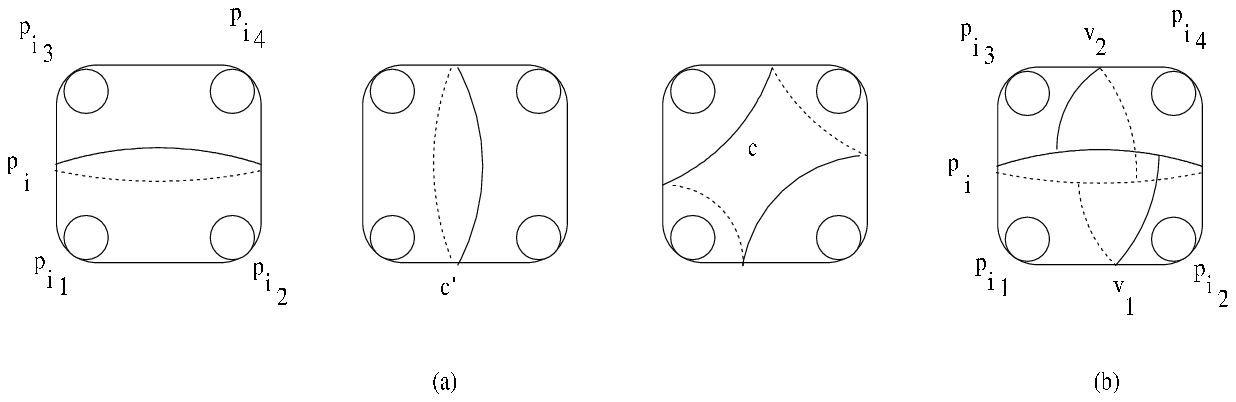}}
\vglue 6pt
\centerline{\small Figure 2.4}  
\endinsert

\section{Fenchel-Nielsen coordinates of Teichm\"uller space}

In this section, we will recall the definition of the Fenchel--Nielsen
coordinates on Teichm\"uller space. The definition below is tailored to 
our
purposes and differs slightly from the usual one (for instance in 
[\ref{IT}]), but
they are equivalent. The basic setup for the Fenchel--Nielsen 
coordinates is
a surface with a colored hexagonal decomposition. The difficulty in 
defining
the coordinates is due to the change in the underlying surfaces as the 
metric
varies in Teichm\"uller space. 

\sh{3.1\qua Marked surfaces}

Recall that a marking on an oriented surface $F$ is colored hexagonal
decomposition $m =(p, b, col)$ of the surface. A marked surface is a 
pair
$(F, m)$ where $m$ is a marking. Two marked surfaces $(F, m)$ and $(F', 
m')$
are \it equivalent \rm if there is an orientation preserving 
homeomorphism
$h\co  F \to F'$ so that $h(m)$ is isotopic to $m'$. It is clear from the
definition that a self-homeomorphism $h\co  F \to F$ is isotopic to the 
identity
if and only if $h(m)$ is isotopic to $m$. A marked hyperbolic surface 
is a
triple $(F, m, d)$ where $(F, m)$ is a marked surface and $d$ is a 
hyperbolic
metric on $F$ with geodesic boundaries. Two marked hyperbolic surfaces
$(F, m, d)$ and $(F', m', d')$ are \it equivalent \rm if there is an
orientation preserving isometry $h\co  F \to F'$ so that $h(m)$ is 
isotopic to
$m'$.

Fix a marked surface $(F, m_{0})$. The Teichm\"uller space of the 
marked
surface, denoted by $T(F)$ is the space of all equivalence classes of 
marked
hyperbolic surface $(G, m, d)$ so that $(G, m)$ is equivalent to $(F, 
m_0)$.

\sh{3.2\qua Metric twisting}

To define the Fenchel--Nielsen coordinate, we will first need the 
following
well known lemma. See [\ref{Bu}] (lemma 1.7.1) for a proof.

\proclaim{Lemma 3.1}Let $F_{0,3}$ be the 3--holed sphere with boundary
components\break $\partial_1, \partial_2, \partial_3$.

\items
\item{\rm(a)}
For any three positive real numbers $x_1, x_2, x_3$, there exists a
hyperbolic metric $d$ on $F_{0,3}$ so that the boundary components
$\partial_i$ are geodesics of lengths $x_i$. Furthermore, the metric 
$d$ is unique up to isometry.

\item{\rm(b)}If the distinct pairs of geodesic boundary components in (a) are 
joined
by the shortest geodesic arcs, then these three arcs are disjoint and 
cut
the surface into two isometric right-angled hexagons.
\enditems\endproc

We also need to introduce the notion of ``metric twisting of a marked
Riemannian annulus along a geodesic'' in order to define the 
coordinate. Let
$A=[-1, 1] \times S^1$ be an oriented annulus with a Riemannian metric 
$d$
so that the curve $\{0\} \times S^1$ is a geodesic. A \it marking \rm 
on $A$
is the homotopy (rel endpoints) class of a path $a\co  [-1, 1] \to A$ so 
that
$a(\pm 1) \in  \{\pm 1\} \times S^1$. Fix a real number $t$. The metric
$t$--twisting of a marked Riemannian annulus $(A, [a],d)$ is a new 
marked
Riemannian annulus $(A', [a'], d')$ defined as follows. First cut the 
annulus
$A$ open along the geodesic $\{0\}\times S^1$ to obtain two annuli
$A_- =[-1,0]\times S^1$ and $A_+ = [0,1]\times S^1$. Let $S_{\pm}^1$ be 
the
geodesic boundary of $A_{\pm}$ corresponding to $\{0\}\times S^1$ and 
let
$\phi\co  S_- \to S_+$ be the isometry so that $A = A_+ \cup _{\phi} A_-$. 
The
circles $S^1_{\pm}$ have the induced orientations from $A_{\pm}$ and
$\phi$ is orientation reversing. Let
$\psi\co  \{Re^{i\theta} | \theta \in \bold R\} \to S^1_{+}$ be an 
orientation
preserving isometry and $\rho\co  S^1_+ \to S^1_+$ be the $ t$--twisting of
$S^1_+$ which sends $x$ to  $\psi( e^{2 \pi i t} \psi^{-1}(x))$. Define 
the
new annuli $A'$  to be $A_+ \cup _{\rho \phi} A_-$. The Riemannian 
metric $d'$
on $A'$ is the gluing metric. To define the marking, let us represent 
the
original marking $[a]$ by a path $a$ so that
$a(0) = a([-1,1]) \cap (\{0\} \times S^1)$. The new path $a'$ on $A'$ 
is given
by $[a|_{[-1,0]}] * [b] * [a|_{[0,1]}]$ where $[x]$ denotes the image 
of $x$
under the quotient map $A_+ \cup A_- \to A'$, $*$ denotes the
multiplication of paths, and $b$ is the geodesic path of length $|t|$ 
in
$S^1_+$ starting from $\rho(\phi(a(0)))$ and ending at $a(0)$ so that
the orientation of $b$ coincides with that of $S^1_+$ if and only if 
$t>0$.
Note that there is a natural identification of the boundary of $A$ and 
$A'$.
For simplicity, we will assume that $\partial A = \partial A'$ under 
this
identification. There exists an orientation preserving homeomorphism
$h\co  A \to A'$ so that $h|_{\partial A} = id$ and $h(a)$ and $a'$ are 
homotopic
rel endpoints. Thus the marked annuli $(A, [a])$ and $(A', [a'])$ are
equivalent. For simplicity, we will denote $(A', [a'], d')$ by
$T_t(A, [a], d)$, $[a'] = T_t([a])$, and $d' =T_t(d)$.

One can also simplify the marking somewhat as follows. It is well known 
that
each path $a\co  [-1, 1] \to [-1, 1]\times S^1$ with
$a(\pm 1) \in \{\pm 1\} \times S^1$ is relative homotopic to an 
embedded arc.
Also relative homotopic embedded arcs are isotopic by isotopies fixing 
the
endpoints. Thus each marking $[a]$ corresponds to a unique isotopy 
class of
proper arc. For this reason, we will usually represent the marking by 
the
isotopy class. 

It follows from the definition that the following holds.

\proclaim{Lemma 3.2}If $t_1, t_2 \in \bold R$,  then
$ T_{t_1}(T_{t_2}(A, [a], d))$ is isometric to\nl $T_{t_1 +t_2}(A, [a], 
d)$ by
an orientation preserving isometry preserving the marking.\endproc

\rk{3.3}
We now recall the Fenchel--Nielsen coordinates on the Teichm\"uller 
space
$T(F)$ of a marked surface $(F, m)$. Let $N = 3g+r -3$. Given a point
$x =(x_1, t_1, x_2, t_3, ...., x_N, t_N,$ $x_{N+1}, ..., x_{N+r})
\in (\bold R_{>0}\times \bold R)^{N} \times \bold R_{>0}$, we will 
describe
the corresponding hyperbolic metric $(FN)^{-1}(x) =[d] \in T(F)$ as 
follows.

Suppose the marking $m$ is $(p, b, col)$ where $p=p_1\cup\dots \cup 
p_{3g+r-3}$
and $p_{3g+r-3+i}$ is the $i$-th boundary component of $F$. Suppose $P$
is a 
component of $F -p_1 \cup \dots \cup p_{3g+r-3}$  
bounded
by $p_i$, $p_k$ and $p_l$ so that the cyclic order
$i \to k \to l \to i$ coincides with the cyclic orientation on the 
boundary
of its red hexagon. Then we denote this component by $P_{ijk}$.
Note that except for the closed surface of genus 2, only one component 
of
the form $P_{ijk}$ or $P_{ikj}$ can exist. 

Now give each 3--holed sphere $P_{ijk}$ a hyperbolic metric so that
 so that (1) the length of $p_r$ is $x_r$ and (2) each arc in
$b \cap P_{ijk}$ is the shortest geodesic arc perpendicular to the 
boundary.
The red hexagon in $P_{ijk}$ is now represented by a right-angled 
hexagon
$H_{ijk}$.

We construct the hyperbolic surface $(FN)^{-1}(x)$ in two steps. 
Let $x' = (x_1, 0,$ $x_2,0, \dots,$ $x_N, 0, x_{N+1}, \dots, x_{N+1})$ be 
the
point having the same $x_i$--th coordinate as $x$ but zero twisting
coordinates. Then the hyperbolic surface in $T(F)$ having
Fenchel--Nielsen coordinates $x'$ is constructed  as follows. Glue $P_{ijk}$ and 
$P_{irs}$
along $p_i$ by an orientation reversing isometry so that it sends the
red interval $p_i \cap H_{ijk}$ to the red interval $p_i \cap H_{irs}$. 
This gluing produces a new hyperbolic surface $(F', d')$ homeomorphic 
to $F$.
The marking $m' =(p'_1\cup \dots\cup p'_{3g+r-3}, b', col')$ on $F'$ 
comes
from the quotient of $\cup p_i$ and $\cup (b \cap P_{ijk})$ and the red
hexagons $H_{ijk}$. By the construction, the marked surfaces $(F, m)$ 
and
$(F', m')$ are equivalent. This gives the point $(FN)^{-1}(x') \in 
T(F)$.

For a general point
$x \in  (\bold R_{>0}\times \bold R)^{3g+r-3} \times \bold R_{>0}$,
the underlying hyperbolic surface $F''$ having $x$ as its 
Fenchel--Nielsen
coordinates is obtained from $F'$ by performing metric $t_i$ twisting 
on
each Riemannian annulus $N(p_i)$ along the geodesic $p_i$. The marking
$m'' =(p'', b'', col'')$ on $F''$ is defined as follows. The 3--holed 
sphere
decomposition of $F''$ corresponds to the quotient of $\cup_i p_i$ in
$\cup P_{ijk}$. To find the hexagonal decomposition, choose the marking
$m'=(p', b', col')$ on $F'$ so that $ b' \cap N(p'_i)$ consists of two 
arcs
$c_{i_1}, c_{i_2}$. Now each isotopy class $[c_{i_r}]$ in the annulus
$N(p'_i)$ is a marking. The new isotopy class of arcs 
$T_{t_i}([c_{i_r}])$
is represented by an embedded arc $c'_{i_r}$ having the same endpoints 
as
that of $c_{i_r}$. We defines $b''$ to be the quotient of
$(b - \cup_{i} int(N(p_i)))\cup (\cup_{i,r} c'_{i_r})$. Define the 
coloring
of the hexagons in $F'' -p''\cup b''$ by the corresponding coloring of 
$F'$.
By the construction, we see that the marked surface $(F'', m'')$ is 
equivalent
to $(F, m)$. This gives the full description of the Fenchel--Nielsen 
coordinate.

The use of the marking is to identify the homotopy classes of loops and
elements in $CS(F)$ on different surfaces. To be more precise, consider 
the
two marked surfaces $(F',m')$ and $(F'',m'')$ constructed above. By the
construction, there is an orientation preserving homeomorphism
$h\co  F' \to F''$ so that $h(m')$ is isotopic to $m''$. This 
homeomorphism
induces a bijection between $CS(F')$ and $CS(F'')$ as follows. If $a'$ 
is a
curve system in $F'$, then the corresponding curve system $a''$ 
homotopic to
$h(a')$ is obtained in the following procedure. Cut $a'$ open along all
$p_i$'s to obtain a collection of geodesic arcs in $P_{ijk}$. Now 
rejoin these
arcs at the ends points in pairs according to the original cutting 
points by
the oriented geodesic arcs in $p_i$ of length $x_i|t_i|$ from the left 
side
endpoints to the right side endpoints along $p_i$. The resulting curve 
system
is $a''$. It follows from the construction that,
$$l_{d''}([a'']) \leq l_{d'}([a']) +
\smash{\sum_{i=1}^{3g+r-3}}\vrule width 0pt depth 10pt x_i|t_i| I([a], p_i). \tag 3.1$$
The basic result about the Fenchel--Nielsen coordinates is that the map
$FN\co$\break  $T(F) \to (\bold R_{>0}\times \bold R)^{3g+r-3} \times \bold 
R_{>0}$ is a
homeomorphism. See for instance [\ref{IT}] chapter 8, or [\ref{Bu}] chapter 6.

\section{Proof of the main theorem}

We prove the main theorem in this section. There are two  
facts about
hyperbolic polygons used in the proof. These two facts will be 
established in
\S5. In subsections 4.1--4.4, we prove the first inequality (1.1). In 
the 
remaining subsections, we establish (1.2). 

To begin the proof, we fix a marking on the surface and let $FN$ and 
$DT$ be
the associated coordinates on the Teichm\"uller space $T(F)$ and the 
space of
curve systems $CS(F)$.

\rk{4.1}To prove inequality (1.1) for all metrics $[d] \in T(F)$ and
$[a], [b] \in CS(F)$, by the remarks in subsection 2.3, it suffices to 
show
$$ |l_d([a]) - l_d([b])| \leq 6 |FN(d)|$$
whenever $DT(a)$ and $DT(b)$ differ only in one intersection coordinate
$x_i$ by 2, ie, $x_i(a) = x_i(b) + 2$ and $x_j(a) = x_j(b)$ for all
$j \neq i$ and $t_k(a) = t_k(b)$ for all $k$. There are three subcases 
we have
to consider according to the nature of the decomposing loop $p_i$:
(1) $[p_i] \in CS(F)$ and is adjacent to only one 3--holed sphere 
$P_{iij}$;
(2) $[p_i] \in CS(F)$ and is adjacent to two different 3--holed spheres
$P_{ii_1i_2}$ and $P_{ii_3i_4}$;
(3) $p_i \subset \partial F$. 

\rk{4.2}In the first case, by proposition 2.3, we can write
$a \cong p_j^{e_1} c^{e_2} b$ where $e_1, e_2 \in \{0, \pm 1, \pm 2\}$ 
and
$c$ is as shown in figure 2.3.

We can write the loop $c \cong p_i^{\pm 1} c'$ where $c'$ has zero 
twisting coordinates as shown in figure 2.3. Let $l(S)$ be the length of the shortest geodesic segment 
in
the 3--holed sphere $P_{iij}$ joining the two boundary components 
corresponding
to $p_i$. Then by the definition of the Fenchel--Nielsen coordinates, we 
have
$l_d([c']) \leq x_i(d)|t_i(d)| + l(S)$. This shows
$$\eqalign{ |l_d([a]) - l_d([b])|& \leq l_d( [p_j^{e_1} c^{e_2}])\cr
& \leq 2 l_d([p_j])  + 2 l_d([c])\cr
& \leq 2 x_j(d) + 2 l_d([p_i^{\pm 1} c'])\cr
& \leq  2x_j(d) + 2 x_i(d) + 2 l_d([c'])\cr
& \leq 2 x_j(d) + 2x_i(d)  + 2 x_i(d)|t_i(d)| +  2l(S).}$$
By proposition 5.1, we can estimate the length $l(S)$ in terms of the red
right-angled hexagon inside $P_{iij}$. Thus we obtain,
$$ l(S) \leq  2/x_i(d) + 2/x_j(d) + x_j(d)/2 + 2\log 2.$$
Combining these together, we obtain
$$\eqalign{  |l_d([a]) - l_d([b])|& \leq 4x_j(d) {+} 2x_i(d) {+} 4/x_j(d) {+} 4/x_i(d) {+}
2 x_i(d)|t_i(d)| {+} 4 \log 2\cr
&\leq 4 |FN(d)|\cr
&\leq 2 |FN(d)| |DT(a) - DT(b)|.}$$

\rk{4.3}In the second case, we use proposition 2.4. Thus
$a \cong p_{i_1}^{s_1}\dots p_{i_4}^{s_4} c^e b$ where $|s_1| + |s_2| 
\leq 2$,
$|s_3| + |s_4| \leq 2$, $e \in \{\pm 1\}$ and $c$ has Dehn--Thurston
coordinates of the form $( [0,0],\dots , [0,0],[2,t],[0,0],\dots, 0)$ 
where
$|t| \leq 2$. See figure 2.4. By the triangle inequality,
$$ |l_d([a]) - l_d([b])| \leq  2\smash{\sum_{j=1}^4}\vrule
width 0pt depth 10pt l_d([p_{i_j}]) + 
l_d([c]).$$
To estimate $c$, let $c' \cong c_{zt}$. Then $c \cong p_i^{ t} c'$ 
where
$|t| \leq 2$ hence $l_d([c]) \leq l_d([c']) + 2 x_i(d).$

Consider the metric $d'$ on $F$ so that $FN(d)$ and $FN(d')$ are the 
same
except at the $i$-th twisting coordinate where $t_i(d') = 0$. Then by 
the
definition of the Fenchel--Nielsen coordinate
$l_d([c']) \leq l_{d'}([c']) + 2 x_i(d) |t_i(d)|$. We will estimate the length
$l_{d'}([c'])$ as follows. Let $v_1$ and $v_2$ be the shortest arcs in 
the
red-hexagons $H_{ii_1i_2}$ and $H_{i i_3 i_4}$ joining the $p_i$--side 
to its
opposite side (see figure 2.4(b)). Then by the construction of the
Fenchel--Nielsen coordinates, we have
$l_{d'}([c']) =  l_{d'}(v_1) +  l_{d'}(v_2)$. By proposition 5.1, we can 
estimate
the lengths $l_{d'}(v_k)$ for $k=1,2$ as follows. For simplicity, we 
write
$x_r = x_r(d)$.
$$\eqalign{ l_{d'}(v_1)& \leq  2/x_i+ 2/x_{i_1} + x_{i_1}/2 + x_{i_2}/2 + \log 
2.\cr
l_{d'}(v_2)& \leq 2/x_i+ 2/x_{i_3} + x_{i_3}/2 + x_{i_4}/2 + \log 
2.}$$
Combining the above formulas, we obtain
$$\eqalign{| l_d&([a]) - l_d([b])|\cr
&\leq 2 \smash{\sum_{j=1}^4}\vrule width 0pt height 15pt  x_{i_j}  + 2x_i + 
x_i(d)|t_i(d)|
+ 4/x_i + 2/x_{i_1} + 2/x_{i_3} + x_{i_1} +\cr
&\kern 3in \smash{x_{i_3} + x_{i_2} + x_{i_4}
+ 4 \log 2}\cr
&\leq  6|FN(d)|\cr
& \leq 3|FN(d)| |DT(a) - DT(b)|.}$$
Note the coefficient is $6$ instead of $4$ since $i_1, i_2, i_3$, and 
$i_4$
need not be distinct indices.

\rk{4.4}In the third case that $x_i(a) = x_i(b) +2$ where
$p_i \subset \partial F$, the result follows from the previous case by 
the
standard metric double construction. Indeed, let $F^*$ be the double of 
$F$
across its boundary, ie, $F^* = F \cup_{id} F$ where $id$ is the 
identity
map on $\partial F$. We give $F^*$ the double metric $d^*$  and the 
marking
the double of the original marking. The double of a curve system
$\alpha \in CS(F)$ is denoted by $\alpha^* \in CS(F^*)$. Note that the
twisting coordinate of $\alpha^*$ at each boundary component is always 
zero.
Then it follows from the definition that $|FN(d^*)|\leq 2|FN(d)|$, and
$|DT([a]^*) - DT([b]^*)| = 2$. Thus by the boundaryless case,
$$|l_{d}([a]) - l_d([b])| = 1/2 |l_{d^*}([a]^*) - l_{d^*}([b]^*)|
\leq 3 |FN(d^*)|$$
$$\leq 6|FN(d)| = 3 |FN(d)| |DT(a) -DT(b)|.$$

\rk{4.5}To prove the second inequality (1.2), we first consider the two 
cases
$FN(d_1) - FN(d_2) = (0,\dots , 0,c,0,...0) \in (\bold R_{>0} \times 
\bold R)^N
\times R_{>0}^r$ where either $c$ is $t_i(d_1) - t_i(d_2)$ or is
$x_j(d_1) - x_j(d_2)$. The general case follows by a simple 
interpolation.
These two cases will be dealt separately.

\rk{4.6}In the first case that $c = t_i(d_1) - t_i(d_2)$, then the metric 
$d_2$
is obtained from $d_1$ by a metric twisting of signed length $x_i(d_1) 
c$.
Thus if $a \in \alpha$ is a $d_1$--geodesic representative, then a
representative $a' \in \alpha$ in the $d_2$--surface is obtained from
$a$ by cutting $a$ open along $p_i$ and gluing $I(\alpha, p_i)$ many 
copies
of geodesic segments of lengths $x_i(d_1) |c|$ as obtained in the 
inequality
(3.1).  Thus
$$| l_{d_1}(\alpha) -l_{d_2}(\alpha)| \leq x_i(d_1) |c| |DT(\alpha)| 
\leq
 D(FN(d_1),FN(d_2)) |DT(\alpha)|.$$

\rk{4.7}In the second case that $c = x_i(d_1) - x_i(d_2)$, due to 
symmetry, it
suffices to show that
$$ l_{d_2} (\alpha) \leq l_{d_1}(\alpha) + 4 D(FN(d_1),FN(d_2)) 
|DT(\alpha)|.$$
To this end, take a $d_1$--geodesic representative $a \in \alpha$. We 
will
construct a piecewise geodesic representative $a' \in \alpha$ in 
$d_2$--surface
and estimate the length $l_{d_2}(a')$. The $d_2$--surface $F'$ is 
obtained from
the $d_1$--surface by cutting open along the geodesic $p_i$. Then 
replace the
3--holed spheres $P_{ijk}$ and $P_{irs}$ adjacent to $p_i$ by new pairs 
so that
the lengths at $p_i$ are $l_{d_2}([p_i])$, and all other lengths remain 
the
same. For each 3--holed sphere $P$ in the decomposition, let
$H$ in $P$ be one of the right-angled hexagon obtained from lemma 
3.1(b). Note
that the metric gluing to obtain the $d_2$--surface has the same 
twisting angles
$t_j$. This shows that there is an orientation preserving homeomorphism 
$h$
from the $d_1$--surface to the $d_2$--surface so that (1) $h$ sends the
right-angled-hexagon $H$ to the right-angled-hexagon $H$; (2) $h$ on 
each edge
in the  boundary of the right-angled hexagons $H$ and $P-H$ are 
homothetic
maps. (Note that the red-hexagons used as part of a marking on the
$d_k$--surface are in general different from the hexagons $H$.) The
representative $a'$ is choosen so that on 
each
right-angled hexagon $X$ = $H$ or $P-H$, $a'$ consists of geodesic segments
and for each component $b$ of $a \cap X$, there exists exactly one 
component
$b'$ of $a' \cap X$ for which $ h(\partial b) = \partial b'$.

It follows from the construction that $l_{d_2}(b') = l_{d_1}(b)$ unless
$b$ lies in either $P_{ijk}$ or $P_{irs}$. In the later case, by 
theorem 5.2,
we have
$$l_{d_2}(b') \leq l_{d_1}(b) + 4 |\log \sinh(x_i(d_1)/2)
- \log \sinh(x_i(d_2)/2)|.$$
Let $n$ be sum of the number of components of $a \cap X$ for all 
right-angled
hexagons $X$ in $P_{ijk}$ and $P_{irs}$. Then
$$ l_{d_2}(\alpha) \leq l_{d_2}(a') \leq l_{d_1}(\alpha) + 4n 
|\log \sinh(x_i(d_1)/2) - \log \sinh(x_i(d_2)/2)|.$$
It remains to estimate the number $n$.

\proclaim{Lemma 4.1}Under the above assumptions
$$ n \leq (|t_i(d_1)| + |t_j(d_1)| + |t_k(d_1)| + |t_r(d_1)| + 
|t_s(d_1)|
+ 7) |DT(\alpha)|.$$\endproc

Assuming this lemma, then we obtain the required estimate that
$$\eqalign{l_{d_2}(\alpha)& \leq l_{d_2}(a')\cr
&\leq l_{d_1}(\alpha) + 4 (|t_i|+ |t_j| + |t_k| + |t_r| + |t_s|
+7) |\log \sinh(x_i(d_1)/2)\cr 
&\kern 2.5in - \log \sinh(x_i(d_2)/2)||DT(\alpha)|\cr
&\leq l_{d_1}(\alpha) + 4 D(FN(d_1),FN(d_2)) |DT(\alpha)|}$$
where $t_n = t_n(d_1)$. Thus the inequality (1.2) follows in this case.

\medskip
{\bf Proof of lemma 4.1}\qua Let us first consider the special case that
$t_j(d_1) =0$ for all $j$. In this case the red-hexagons in the 
$d_1$--surface
are the same as the right-angled hexagon $H$. Thus
$n \leq I(\alpha, p) + I(\alpha, b)$ where $(p, b, col)$ is the marking 
on
the $d_1$--surface. Now we can write
$\alpha = [p_1^{r_1}.... p_N^{r_N}]\alpha_{zt}$ where $r_i$ is the
Dehn--Thurston twisting coordinate of $\alpha$ and $\alpha_{zt}$ has 
zero
twisting coordinates. Thus,
$$\eqalign{ n &\leq I(\alpha, p) + I(\alpha_{zt}, b) + I(p_1^{|r_1|} ... 
p_N^{|r_N|}, b)\cr
&\leq 2I(\alpha, p) + 2 \sum_{i=1}^N |r_i|\cr
&\leq 2 |DT(\alpha)|.}$$
In particular, the conclusion holds in this case. Also we see that for 
any
marking $(p, b, col)$ on a surface,
$I(\alpha, p) + I(\alpha, b) \leq 2 |DT(\alpha)|$.

In the general case that some $t_j(d_1) \neq 0$, we take all $p_j$'s to 
be
$d_1$--geodesics
and let $u_{hl}$ be the shortest geodesic segment joining $p_h$ to 
$p_l$ when
$p_h$ and $p_l$ lie inside some 3--holed sphere component of $F-p$. Let
$b$ be the $d_1$--geodesic representative of the marking curve and 
$b_{hl}$ be
the component of $b \cap P_{hlm}$ corresponding to $u_{hl}$. Then by
definition of Fenchel--Nielson coordinates, $u_{hl}$ is relatively 
homotopic
to $ w_h * b_{hl} * w_l$ where $w_h$ is a geodesic path in $p_h$ of 
length $
x_h(d_1) |t_h(d_1)|$. Thus the number of new intersection points in
$a \cap w_h$ is at most $(|t_h(d_1)| + 1) I(\alpha, p_h)$. This shows 
that
$$\eqalign{ n &\leq | a \cap p| + \sum_{h,l} |a \cap (\cup_{h,l} u_{hl})|\cr
&\leq I(\alpha, p) + |a \cap b| + \sum_{h} (|t_h(d_1)|+1) I(\alpha, 
p_h)\cr
&\leq I(\alpha, p) + I(\alpha, b) + \sum_{h} (|t_h(d_1)|+1) I(\alpha, 
p)\cr
&\leq ( \sum_{h} |t_h(d_1)| + 7) |DT(\alpha)|,}$$
where the sum is over the set
$\{i,j,k,r,s\}$. 

\rk{4.8}The above estimate works even if the loop $p_i$ is a boundary 
component
of the surface $F$.

\sh{4.9\qua The general case} 

The general case of any two metrics $d_1$ and $d_2$ follows from
 interpolation. Namely we use the formula
$|F(x_i, t_i) - F(y_i, s_i)| \leq |F(x_i, t_i) - F(y_i, t_i)| +
| F(y_i, t_i) - F(y_i, s_i)|$. Thus the result follows. Also the 
corollary
1.2 follows from the standard argument involving the definition of the 
space of measured laminations. See [\ref{LS}] \S6 for the proof of the similar result 
for the intersection pairing.\endproof

\section{Elementary facts about hyperbolic polygons}

We will prove two facts used in the proof of the main theorem in this
section. For basic information on hyperbolic hexagons, see [\ref{Be}] \S7.19, [\ref{Bu}] \S2.4.
Suppose $H$ is a right-angled hyperbolic hexagon whose side 
lengths (reading from counterclockwise) are : $a, z, x, y, b$ and $w$. See 
figure 1.1.

\proclaim{Proposition 5.1}Consider the right-angled hexagon $H$ above.
Let $h$ be the length of the shortest geodesic arc from the $a$--side to 
the
$y$--side. Then:

\items\item{\rm(a)}$w \leq 1/a + 1/b + x + 2\log 2.$

\item{\rm(b)} $ h \leq 1/a + 1/b + b + x + \log 2$ and
$ h \leq  1/a + 1/2( 1/b+ 1/x) + b + x + 2\log 2.$\enditems\endproc

\prf By the cosine rule,
$\cosh w = (\cosh x +\cosh a \cosh b)/(\sinh a \sinh b)$.
Using $\cosh a \cosh b + \cosh c \leq \cosh a \cosh b(\cosh c + 1)$ and
$\cosh w \geq 1/2 e^{w}$, we obtain
$$ 1/2 e^{w} \leq \coth a \coth b (\cosh x + 1).$$
Taking logs, we get
$$ w - \log 2 \leq \log \coth a + \log \coth b + \log ( \cosh x + 1).$$
On the other hand, $\coth a \leq 1+ 1/a$. Thus
$\log \coth a \leq \log (1+ 1/a) \leq 1/a$.  Similarly, $\log \coth b 
\leq 1/b$.
Finally, $\log( \cosh x + 1) \leq \log ( e^x +1) \leq x + \log 2$. Put 
all
these together, we obtain the estimate (a).

To see (b), by the cosine law for pentagon,
$$ \cosh h = \sinh b \sinh w.$$
Now $ e^h/2 \leq \cosh h$ and $\sinh x \leq e^x/2$. Thus
$e^h \leq 1/2 e^b e^{w}$. This shows that $h \leq b + w -\log 2$. By 
part
(a), we obtain
$$ h \leq 1/a + 1/b + x + b + 2 \log 2.$$
Also $h \leq 1/a + 1/x + x + b + \log 2$. Thus
$$h \leq 1/a + 1/2(1/b + 1/x) + b + x + \log 2.\eqno{\sq}$$

Let $A_g = A_g(\lambda, \mu)$ be a geodesic segment in $H$ joining two 
sides
of $H$ so that the endpoints of $A_g$ cut the sides into two intervals 
of
lengths $\lambda t$, $(1-\lambda)t$ and $\mu r$, $(1-\mu) r$. In the
discussion below, the numbers $a,b,\lambda, \mu$ remain constant. The
variable is $x$ and $y,z, w$ depend on $x$. Let $S = S_{\lambda, \mu}$ 
be the
length of $A_g$. Our goal is to estimate the rate of change of
$S_{\lambda, \mu}$ with respect to $x$.

\proclaim{Theorem 5.2} 
Under the above assumption, we have $|\frac{dS}{dx}| \leq 4 \coth x$. 
\endproc

\prf We begin with several simple lemmas based on the cosine 
and
sine laws in hyperbolic geometry.

\proclaim{Lemma 5.3}In the right-angled hexagon $H$,

\items
\item{\rm(a)}  $\frac{dy}{dx} = -\frac{ \coth z}{\sinh x} $.

\item{\rm(b)} $\frac{dw}{dx} = \frac{1}{\sinh a \sinh z}$.

\item{\rm(c)} $ | \frac{dy}{dx}| < \coth  x$ and $\frac{dy}{dx} < 0$.

\item{\rm(d)}  $ 0 < \frac{dw}{dx} < 1$.

\item{\rm(e)} $\cosh x > \coth z$  and $\sinh x \sinh z > 1$.
\enditems\endproc

\prf From the cosine rule:
$\cosh x = (\cosh w + \cosh y \cosh z)/(\sinh y \sinh z)
       = \cosh w/(\sinh y \sinh z) + \coth y \coth z$
$
       > \coth y \coth z   
       > \coth z.   $
Now squaring the inequality and using $\cosh^2 x = 1 + \sinh^2x$ and
$\coth^2 z = 1 + 1/\sinh^2 z$, we obtain $\sinh x \sinh z > 1$.
This shows $(e)$.

Differentiating the other cosine rule
$$\cosh y = (\cosh a + \cosh b \cosh x)/(\sinh b \sinh x)$$
gives
$\frac{dy}{dx} = - (\cosh b + \cosh a \cosh x)/(\sinh b \sinh^2 x \sinh 
y).$
Plugging in the cosine rule  $\cosh b + \cosh a \cosh x = \sinh a \sinh 
x 
\cosh z$, we obtain,
$$\eqalign{\frac{dy}{dx}& = - (\sinh a \sinh x \cosh z)/
(\sinh b \sinh^2 x \sinh y)\cr
     & = - (\sinh a \cosh z)/(\sinh b \sinh x \sinh y).}$$
Plugging in the sine rule  $\sinh a/\sinh y = \sinh b/\sinh z$  gives
$$\frac{dy}{dx} = - \cosh z /(\sinh x \sinh z) = - \coth z/\sinh x.$$
This shows (a) and the second part of (c).

By the inequality $(e)$ above
$ |\frac{dy}{dx}| < \coth x/\coth y < \coth x.$
This shows (c).

For $\frac{dw}{dx}$,
we have
$\cosh w = (\cosh x + \cosh a \cosh b)/(\sinh a \sinh b).$
By the sine law,
$$\frac{dw}{dx} = \sinh x/(\sinh a \sinh b \sinh w)=1/(\sinh b \sinh y)
=1/(\sinh a \sinh z).$$

By the rewritten form of (e) for the pair $(a,z)$ instead of $(z,x)$, 
we have
 $\sinh a \sinh z$$ > 1$. This shows $ 0 < \frac{dw}{dx} <1. $ Thus  
both
(b) and (d) hold.\hfill$\sq$

\medskip

The next lemma is well known. It is a simple application of the sine
law. We will omit the details of the proof.

\proclaim{Lemma 5.4} Suppose $\Delta qpr$ is a hyperbolic triangle with 
angle
at $p$ being $\alpha$. Suppose starting at time $t=0$ the endpoint $p$ 
moves
along the ray $pr$ with unit speed while the other two points $q, r$ 
remain
fixed. Let $l_{pq}$ denote the length between $p$ and $q$. Then
$dl_{pq}/dt|_{t=0} = - \cos \alpha.$ \qed \endproc

\midinsert
\vglue-15pt\centerline{\epsfbox{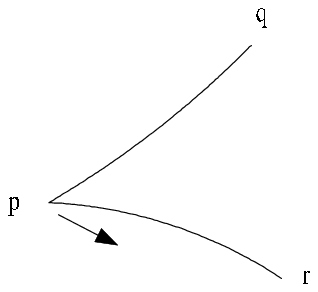}}
\vglue 6pt
\centerline{\small Figure 5.1}  
\endinsert

The next lemma is crucial for most of the estimates in the proof of 
theorem 5.2.

\proclaim{Lemma 5.5} 
Consider a hyperbolic quadrilateral with side lengths and angles 
(reading from
counterclockwise) as $ c$ (side),  right angle, $t$ (side), 
right-angle, $e$
(side), $\beta$ (angle), $S$ (side) and $\alpha$ (angle). 
Consider varying  $t$ and holding $c$ and $e$  fixed, then  

\centerline{ $ 0 < \frac{\partial S}{\partial t} < \coth(t/2).$}
\endproc

\prf By the cosine law,
$\cosh S = - \sinh c \sinh e + \cosh c \cosh e \cosh t.$
Differentiating this equation gives
$\partial S/\partial t = \cosh c \cosh e \sinh t/\sinh S > 0.$
Plugging in the identity  $\cosh v = 2 \sinh^2(v/2) + 1$  three times 
to the
above cosine law gives
$\sinh^2(S/2) = \sinh^2((c-e)/2) + \cosh c \cosh e \sinh^2(t/2) 
            > \cosh c \cosh e \sinh^2(t/2).$
Using $ \sinh t =$$ 2 \coth(t/2)$$ \sinh^2(t/2)$, we obtain the result.
\endproof

We now begin the proof of the theorem 5.2. We will break it into three 
cases,
each of which will have several subcases. We refer to the case where
the geodesic segment $A_g$ has endpoints on adjacent sides as case 1,
sides two apart as case 2 and endpoints on opposite sides as case 3.
In the following discussion, we will assume the hexagon has side 
lengths
$x, y(x), b, w(x), a, z(x) $ where  a  and  b  are fixed. We will use
$\frac{dy}{dx} $ etc, for derivatives of these side lengths. When 
looking
at $S(x)$  however, we will often consider $S$  as a side of a 
hyperbolic
polygon with the angles not incident on  $S$  all right angles. In such 
a
case, we can vary the other sides independently and we will use
$\partial S/\partial c $ for the change in $S$ when we vary only the 
side $c$
of this polygon. 

{\bf Case 1}\qua There are up to symmetry three subcases depending on which
sides $S$ joins, however we will do all three cases simultaneously with 
a
little care.  In this case consider the right-angled triangle cut out 
by
the segment $A_g(\lambda, \mu)$. The side lengths of the triangle  are
$\mu c, \lambda e $ and  $ S$, where $(c, e)$ may be $(x,y), (y,b)$ or
$(b,w)$. Let $\alpha $ be the angle opposite $ \mu c $ and $ \beta $ 
the
angle opposite $\lambda e$. By lemma 5.4, if $c$ increases one endpoint 
$S$
moves off at an angle of $\pi - \beta,$ hence
$\partial S/\partial c = \mu \cos(\beta)$, similarly as $e$ increases 
the
other endpoint of $S$ moves off at an angle of $\pi-\alpha,$ hence
$\partial S/\partial e = \lambda \cos(\alpha).$ Thus
$dS/dx = (\partial S/\partial c)(dc/dx) + (\partial S)/(\partial e) 
(de/dx)
      = \mu \cos(\beta) (dc/dx) + \lambda \cos(\alpha) (de/dx).$
Since $ 0< \alpha, \beta  < \pi/2$, the cosines are positive. In any of 
the
three cases for  $(c,e)$, by lemma 5.3,  we have $ (dc/dx)(de/dx)\leq 0 
$.
Therefore, by lemma 5.3 again,
$$\eqalign{|dS/dx| &\leq  \max( \mu \cos(\beta) |dc/dx|, \lambda \cos(\alpha) 
|de/dx|)\cr
        &\leq  \max( |dc/dx|, |de/dx|)
        < \coth x.}$$

\midinsert
\centerline{\epsfxsize 4.2in\epsfbox{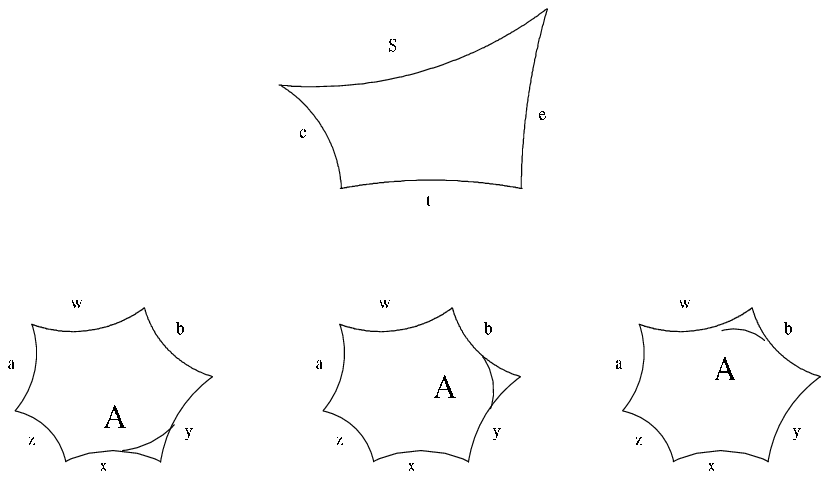}}
\vglue 6pt
\centerline{\small Figure 5.2}  
\endinsert

{\bf Case 2}\qua This case splits into four subcases up to symmetry. We will at
least start these cases together. We have a quadrilateral with sides 
and angles (reading from counterclockwise) as $ \mu c$ (side), 
right-angle,
$t$ (side), right-angle, $\lambda e$ (side), $\beta$ (angle), $S$ 
(side), and
$\alpha$ (angle). Here $(c,t,e)$ is one of $(z,x,y)$, $(x,y,b),$ 
$(y,b,w),$
or $(b,w,a)$.

By Lemma 5.4, $\partial S/\partial c = \mu \cos(\alpha)$  and 
$\partial S/\partial e = \lambda \cos(\beta)$. Note that both of these 
have
magnitude at most  1. Combining this fact with Lemma 5.5, we obtain
$$|dS/dx| = | (\partial S/\partial c)( d c/ dx )
+ (\partial S/\partial e)( de /dx)
+ (\partial S/\partial t)( dt/dx)| $$
$$ \leq |dc/dx| + |de/dx| + \coth(t/2) |dt/dx|.$$
In any case, by lemma 5.3,  $|dc/dx| < \coth x$  and  $|de/dx| < \coth 
x$.
Hence $|dS/dx| \leq  2 \coth x + \coth(t/2) |dt/dx|.$

Subcase (i). $(c,t,e) = (z,x,y)$. In this case $ t=x, dt/dx=1$  and 
using 
the fact that  $2 \coth x = \coth(x/2) + \tanh(x/2) > \coth(x/2) $ we 
see
$|dS/dx| < 4 \coth x.$

Subcase (ii). $(c,t,e) = (x,y,b)$. In this case $ t=y,$  $de/dx=0$,  
and by
lemma 5.3 we have
$$\eqalign{|dS/dx| &\leq |dc/dx| + \coth(t/2) |dt/dx|\cr & <  \coth x + \coth(y/2) 
\coth z/\sinh x
\cr &< \coth x + \coth(y/2) \coth x /\coth y < 3 \coth x.}$$
Note that
$\coth z/\sinh x < \coth x /\coth y$ by the proof of lemma 5.3.

Subcase (iii). $(c,t,e) = (y,b,w)$. In this case  $t=b$  and  
$dt/dx=0.$

Subcase (iv).  $(c,t,e) = (b,w,a).$ In this case  $dc/dx=de/dx=0, t=w $ 
and  
$dt/dx = 1/(\sinh a \sinh z).$ Hence
$$\eqalign{0 &< dS/dx < \coth(w/2)/(\sinh a \sinh z)\cr
&= (1 + \cosh w)/(\sinh w \sinh a \sinh z)\cr
&= (1 + \cosh w)/(\sinh z \sinh x \sinh y)\cr
&< 2 \cosh w/(\sinh z \sinh x \sinh y).}$$
Since 
$\cosh w=\cosh x \sinh y \sinh z - \cosh y \cosh z < \cosh x \sinh y 
\sinh z$,
it follows that $0 < dS/dx < 2 \coth x.$ This completes Case 2.

\midinsert
\centerline{\epsfbox{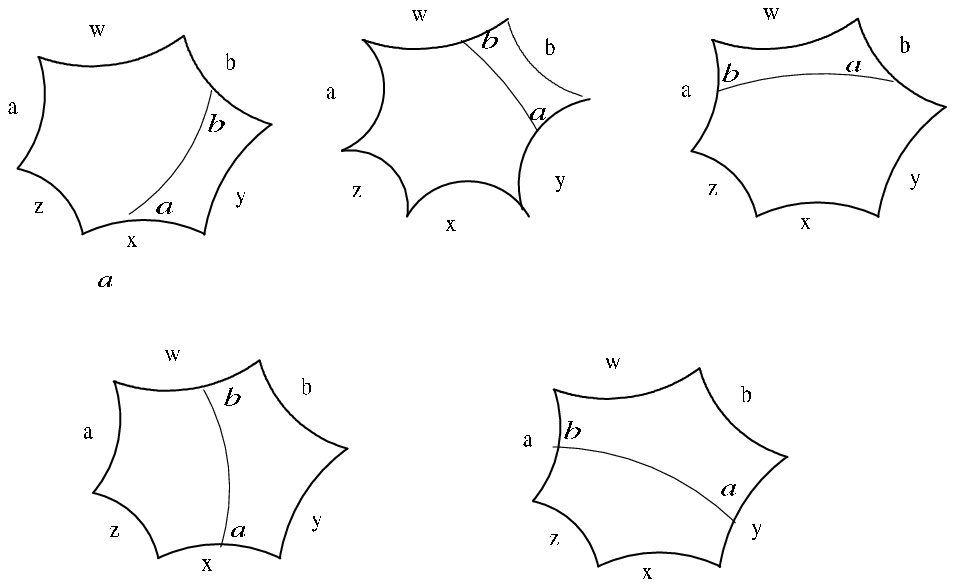}}
\vglue 6pt
\centerline{\small Figure 5.3}
\endinsert

{\bf Case 3}\qua Here there are two subcases (up to symmetry). Either $S$ joins 
$x$ to
$w$ or $S$ joins $a$ to $y$. In the first subcase we have a pentagon 
with
sides and angles (reading from counterclockwise):
$\mu x $ (side), right angle, $y$ (side), right-angle, $b$ (side), 
right-angle,
$ \lambda w $ (side), $\beta$ (angle), $S$ (side) and $\alpha$ (angle).

By Lemma 5.4, $\partial S/\partial x = \mu \cos(\alpha)$ hence
$|\partial S/\partial x| \leq 1$,
and similarly $(\partial S/\partial w) = \lambda \cos(\beta)$  hence
$|\partial S/\partial w| \leq 1.$
Also from Lemma 5.4, increasing $y$  is equivalent to pulling the 
endpoint
of $S$ off at an angle of $(\pi/2)+\alpha $ but $ \cosh(\mu x)$ times 
as
fast, hence 
$(\partial S/\partial y) = \cosh(\mu x) \sin(\alpha) \leq \cosh(\mu 
x).$
Combining these and lemma 5.3, we obtain,
$|dS/dx|$$=$$ |(\partial S/\partial x) $ $+(\partial S/\partial w)$$ 
\frac{dw}{dx}
+ (\partial S/\partial y) \frac{dy}{dx}|$
$\leq  2 + \cosh(\mu x) \coth z/\sinh x.$
To estimate the size, we note that this case is symmetric. On the other 
side 
of $S$ is another pentagon and the same argument gives
$|dS/dx| \leq  2 + \cosh((1-\mu) x) \coth y/\sinh x$. Combining these 
gives
$$|dS/dx| \leq 2 + \min[ \cosh(\mu x) \coth z,
\cosh((1-\mu) x) \coth y]/\sinh x.$$
Since the  min  is at most the geometric mean we get
$$|dS/dx| \leq 2 + [\cosh(\mu x) \cosh((1-\mu) x) \coth z \coth 
y]^{1/2}
/\sinh x.$$
By lemma 5.3, $ \coth y \coth z < \cosh x $ and
$\cosh(\mu x) \cosh((1-\mu) x) = [\cosh(x) + \cosh((1-2\mu)x)]/2 < 
\cosh x.$
Hence we get
$$|dS/dx| < 2 + \cosh x/\sinh x < 3 \coth x.$$

In the second subcase we have two pentagons. One with sides and angles:
$\mu y$ (side), right-angle, $b$ (side), right-angle, $w$ (side), 
right-angle,
$\lambda a $ (side), $\beta$ (angle), $S$ (side) and $\alpha$ (angle).
The other pentagon has  sides and angles: $(1-\lambda) a $(side), 
right-angle,
$z$ (side), right-angle, $x$ (side), right-angle, $  (1-\mu)y $ (side),
$\pi-\alpha $ (angle), $S$ (side), $\pi-\beta$ (angle).

Looking at the first pentagon, by Lemma 5.4, $\partial S/\partial y 
=\mu 
cos(\alpha) $ which has magnitude at most 1. Increasing $w$ by an
infinitesimal amount $\delta w $ has the effect of moving an endpoint 
of $S$
a distance $\cosh(\lambda a)\delta w$ at an angle of $\pi/2 + \beta$. 
Hence 
$\partial S/\partial w = - \cos(\pi/2 + \beta) \cosh(\lambda a).$ and
$dS/dx=\sin(\beta) \cosh(\lambda a)\frac{dw}{dx}+\mu 
\cos(\alpha)\frac{dy}{dx}.$
Note that the first term in always positive and the second may be 
either
positive or negative. Hence we see that
$$dS/dx \geq \frac{dy}{dx} \geq - \coth x.$$
Thus we need only give an upper bound on $ dS/dx.$ The bound above 
gives
$dS/dx \leq \cosh(\lambda a) (\frac{dw}{dx}) + \coth x
      = \cosh(\lambda a)/(\sinh a \sinh z) + \coth x.$
We will derive two upper bounds from this. First since
$$\eqalign{\sinh a &= \cosh(\lambda a) \sinh((1-\lambda) a) 
+ \sinh(\lambda a) \cosh((1-\lambda) a)\cr
&> \cosh(\lambda a) \sinh((1-\lambda) a),}$$
we have
$$ dS/dx < 1/(\sinh((1-\lambda) a) \sinh z) + \coth x  \tag 1$$
Second, since  $\cosh(\lambda a) \leq \cosh a $ and from lemma 5.3 
above we
have $\cosh z > \coth a$, we conclude that
$$dS/dx < \coth z + \coth x  \tag 2$$
Now we turn to the second pentagon to get a third inequality. By Lemmas 
5.4
and 5.5, we see that 
$\partial S/\partial y = (1-\mu) \cos(\pi - \alpha) = -(1-\mu) 
\cos(\alpha)$
and $\partial S/\partial z = \sin(\beta) \cosh((1-\lambda) a)$. Thus
$dS/dx = - (1-\mu) \cos(\alpha) \frac{dy}{dx} + (\partial S/\partial x)
+ \sin(\beta) \cosh((1-\lambda) a) dz/dx.$
Since $dz/dx < 0$, the third term is negative and the first term is at
most $|\frac{dy}{dx}| < \coth x.$ Hence
$$dS/dx < (\partial S/\partial x) + \coth x. \tag 3$$
To bound the first term we want to use Lemma 5.5 above. Let $P$ be the
vertex between $ S$  and  $(1-\lambda) a$. Draw the perpendicular from 
$P$
to $ x$  and call the foot of the perpendicular $ Q$.  Let  $r$ be the
distance from side $z$  to $Q$. Clearly
$  (1-\lambda) > r$ since $r$ is the shortest distance between two 
geodesics.
Applying Lemma 5.5 to the quadrilateral with sides  $PQ, x-r, (1-\mu) y 
$ and
$S$ shows $\partial S/\partial (x-r)< \coth((x-r)/2)$. But
$\partial S /\partial x = \partial S / \partial (x-r)$ as we can make
the infinitesimal change of $(x-r)$ at the end point other than $Q$.
Hence,
$$ dS/dx < \coth((x-r)/2) + \coth x  \tag 4$$

Now we show that if there is $x$ so that $dS/dx> M + \coth x$ for some 
constant
$M$, then $M <3$. Thus $dS/dx \leq 3 + \coth x < 4\coth x$. By (2) we 
see 
$\coth z > M.$ From lemma 5.3, we have $\cosh x >\coth z$. Hence 
$x > arc\cosh M$. From (4) we see  $\coth((x-r)/2) > M $ and hence  
$r > x - 2 arc\coth M.$ Hence $ (1-\lambda) a > r > x - 2 arc\coth M.$
From (1) we have $ 1 > M \sinh((1-\lambda) a) \sinh z$. By lemma 5.3 
again,
we have $\sinh x \sinh z > 1$.  Hence 
$$1 > M \sinh((1-\lambda) a)/\sinh x > M \sinh(x - 2 arc\coth M)/\sinh 
x.$$
Since  $d(\log(\sinh t))/dt = \coth t $ is a decreasing function of  
$t$, we
know  $\sinh(t - c)/\sinh t$  is an increasing function of $t$ therefore
$$1 > M \sinh(arc\cosh(M) - 2 arc\coth M)/\sqrt{M^2 - 1}$$
$$= M(M^2 + 1)/(M^2 - 1) - 2 M^3/ (M^2 - 1)^{3/2}.$$
Thus we get a contradiction if  $M \geq 3$. Thus $dS/dx < 4\coth x$
and we are done.\endproof

\references
\bye